\documentclass[12pt,reqno]{amsart}
\usepackage{amsfonts}

\usepackage{fancyhdr}
\usepackage[latin1]{inputenc}
\usepackage{amsmath}
\usepackage{thmdefs}
\usepackage{natbib}
\usepackage{placeins}
\usepackage{graphicx}
\usepackage[colorlinks=true,linkcolor=red,citecolor=blue]{hyperref}
\hypersetup{
pdftitle={General bootstrap for dual divergences},
pdfauthor={\textcopyright\ Bouzebda and Cherfi},
}
\usepackage{dsfont}
\usepackage{mathrsfs}

\theoremstyle{definition} \theoremstyle{remark}
\numberwithin{equation}{section}

\renewcommand{\cite}{\citet}
\numberwithin{equation}{section}

\begin{document}

\title{General bootstrap for dual $\phi$-divergence estimates}
\author{Salim Bouzebda$^{1,2}$ 
 and Mohamed Cherfi$^{2}$
}
\address{$^{1}$
Laboratoire de math\'ematiques appliqu\'ees \\universit\'e de technologie de Compi\`egne, B.P. 529\\ 60205 Compi\`egne cedex, France}
\address{
$^{2}$L.S.T.A.,  Universit\'e Pierre et Marie Curie\\4 place Jussieu
    75252 Paris Cedex 05, France}
    
    \address{ E-mail adresses : salim.bouzebda@upmc.fr ; mohamed.cherfi@gmail.com}

\date{\today}
\maketitle
\begin{abstract} A general notion of bootstrapped $\phi$-divergence estimates constructed by exchangeably
 weighting sample is introduced. Asymptotic properties of these generalized bootstrapped $\phi$-divergence
 estimates are obtained, by mean of  the empirical process theory,  which are applied to construct the bootstrap confidence set with asymptotically correct coverage probability.
  Some of practical problems are discussed, including in particular, the choice of escort parameter  and several examples of divergences are investigated. Simulation results are provided to illustrate the finite sample performance of the proposed estimators. 
\end{abstract}
\section{Introduction}
\noindent The $\phi$-divergence modeling has proved to be a  flexible tool and provided a powerful statistical
 modeling framework in a variety of applied and theoretical contexts [refer to 
 \cite{BroniatowskiKeziou2009}, \cite{Pardo2006} and  \cite{LieseVajda2006,LieseVajda1987} and the references therein].
 For good recent sources of references to research literature in this area along with statistical applications consult 
 \cite{BasuShioyaPark2011} and \cite{Pardo2006}. Unfortunately, in general, the limiting distribution  of the estimators, or their functionals,
 based on $\phi$-divergences depend crucially on the unknown distribution, which is a serious
  problem in practice. To circumvent this matter, we shall propose, in this work, a general bootstrap
  of $\phi$-divergence based estimators and study some of its properties by mean of a sophisticated
   empirical process techniques.\\
A major application for an estimator is in the calculation of
confidence intervals. By far the most favored confidence interval
is the standard confidence interval based on a normal or a Student
$t$-distribution. Such standard intervals are useful tools, but
they are based on an approximation that can be quite inaccurate in
practice. Bootstrap procedures are an attractive alternative. One
way to look at them is as procedures for handling data when one is
not willing to make assumptions about the parameters of the
populations from which one sampled. The most that one is willing
to assume is that the data are a reasonable representation of the
population from which they come. One then resamples from the data
and draws inferences about the corresponding population and its
parameters. The resulting confidence intervals have received the
most theoretical study of any topic in the bootstrap analysis.\\
Our main findings, which are analogous to that of
\cite{ChengHuang2010},  are summarized as follows. The
$\phi$-divergence estimator
$\widehat{\boldsymbol{\alpha}}_{\phi}(\boldsymbol{\theta})$ and
the bootstrap $\phi$-divergence estimator
$\widehat{\boldsymbol{\alpha}}^{\ast}_{\phi}(\boldsymbol{\theta})$
are obtained by optimizing the objective function
$h(\boldsymbol{\theta},\boldsymbol{\alpha})$ based on the
independent and identically distributed [i.i.d.] observations
$\mathbf{ X}_{1},\ldots,\mathbf{ X}_{n}$ and the bootstrap sample
$\mathbf{ X}_{1}^{\ast},\ldots,\mathbf{ X}_{n}^{\ast}$, respectively,
\begin{eqnarray}
\widehat{\boldsymbol{\alpha}}_{\phi}(\boldsymbol{\theta}):=\arg\sup_{\boldsymbol{\alpha}\in
 \mathbf{\Theta}}\frac{1}{n}\sum_{i=1}^{n}h(\boldsymbol{\theta},\boldsymbol{\alpha},\mathbf{ X}_{i}),\\
\widehat{\boldsymbol{\alpha}}^{\ast}_{\phi}(\boldsymbol{\theta}):=\arg\sup_{\boldsymbol{\alpha}
\in \mathbf{\Theta}}\frac{1}{n}\sum_{i=1}^{n}h(\boldsymbol{\theta},\boldsymbol{\alpha},\mathbf{ X}_{i}^{\ast}),
\end{eqnarray}
where $\mathbf{ X}_{1}^{\ast},\ldots,\mathbf{ X}_{n}^{\ast}$ are independent draws with replacement from the original
 sample. We shall mention that $\widehat{\boldsymbol{\alpha}}^{\ast}_{\phi}(\boldsymbol{\theta})$ can alternatively be expressed as 
\begin{eqnarray}
\widehat{\boldsymbol{\alpha}}^{\ast}_{\phi}(\boldsymbol{\theta})=\arg\sup_{\boldsymbol{\alpha}\in \mathbf{\Theta}}\frac{1}{n}\sum_{i=1}^{n}W_{ni}h(\boldsymbol{\theta},\boldsymbol{\alpha},\mathbf{ X}_{i}),
\end{eqnarray}
where the bootstrap weights are given by $$(W_{n1},\ldots,W_{nn})\sim ~~\mbox{ Multinomial}(n; n^{-1},\ldots,n^{-1}).$$
In this paper, we shall consider the more general exchangeable bootstrap weighting
scheme that includes Efron's bootstrap [\cite{Efron79} and \cite{EfronTibshirani1993}]. The general
resampling scheme was first proposed in \cite{Rubin1981} and extensively studied by
\cite{BickelFreedman1981}, who suggested the name
``weighted bootstrap'', e.g.,
Bayesian Bootstrap when $(W_{n1},\ldots,W_{nn})=(D_{n1},\ldots,D_{nn})$ is equal in
 distribution to the vector of $n$ spacings of $n-1$ ordered uniform $(0,1)$ random variables,
 that is
$$
(D_{n1},\ldots,D_{nn})\sim ~~\mbox{ Dirichlet}(n;1,\ldots,1).$$
The interested reader may refer to \cite{Lo1993}.
The case $$
(D_{n1},\ldots,D_{nn})\sim ~~\mbox{ Dirichlet}(n;4,\ldots,4)$$ was considered in
\cite[Remark 2.3]{Weng1989} and \cite[Remrak 5]{Zheng1988}. The Bickel and Freedman result
concerning the empirical process has been subsequently generalized for empirical processes based on observations in
 $\mathbb{R}^d$, $d>1$ as well as in very general sample spaces and for various set and function-indexed random objects
  [see, for example \cite{Beran1984}, \cite{BeranMillar1986}, \cite{BeranLeCamMillar1987}, \cite{Gaenssler1992}, \cite{Lohse1987}].
    In this framework, \cite{CsorgoMason1989} developed similar results for a variety of other statistical functions. This line of
     research was continued  in the work of \cite{GineZinn1989,GineZinn1990}.  There
is a huge literature on the application of the bootstrap
methodology to nonparametric kernel density and regression
estimation, among other statistical procedures, and it is not the
purpose of this paper to survey this extensive literature. This
being said, it is worthwhile mentioning that the bootstrap as per
Efron's original formulation (see \cite{Efron79}) presents some
drawbacks. Namely, some observations may be used more than once
while others are not sampled at all. To overcome this difficulty,
a more general formulation of the bootstrap has been devised: the
{\it weighted} (or {\it smooth}) bootstrap, which has also been
shown to be computationally more efficient in several
applications. We may refer to  \cite{MasonNewton1992}, \cite{PraestgaardWellner1993} and \cite{delBarrioMatran2000}. \cite{HolmesReinert2004}
 provided new proofs for many known results about the convergence
in law of the bootstrap distribution to the true distribution of smooth
statistics employing the techniques based on Stein's method for empirical processes. Note that other variations of
 Efron's
bootstrap are  studied in \cite{ChatterjeeBose2005} using the term ``generalized bootstrap''.
The practical usefulness of the more general scheme is well-documented in
the literature. For a survey of further results on weighted
bootstrap the reader is referred to \cite{Bertail95}.

The remainder of this paper is organized as follows. In the forthcoming section we recall the
estimation procedure based on $\phi$-divergences. The bootstrap of $\phi$-divergence estimators are introduced, in details, and 
their asymptotic properties are given in Section \ref{sec2}. In Section \ref{Examples}, we provide some examples explaining the computation of the  $\phi$-divergence estimators. In Section \ref{censored},
we illustrate how to apply our results in the context of right censoring. 
Section \ref{simulation} provides simulation results in order to illustrate the performance of the proposed
estimators. To avoid
 interrupting the flow of the presentation, all mathematical developments are relegated to the Appendix.

\section{Dual divergence based estimates}

\noindent The class of \emph{dual} divergence estimators has been recently introduced by
 \cite{Keziou2003} and \cite{BroniatowskiKeziou2009}. 
 Recall that the $\phi$-divergence between a bounded
signed measure $\mathbb{Q}$, and a probability measure $\mathbb{P}$ on
$\mathscr{D}$, when $\mathbb{Q}$ is absolutely continuous with
respect to $\mathbb{P}$, is defined by
$$D_\phi(\mathbb{Q},\mathbb{P}):=\int_{\mathscr{D}}
\phi\left(\frac{{\rm d}\mathbb{Q}}{{\rm d}\mathbb{P}}\right)~{\rm d}\mathbb{P},$$
where $\phi(\cdot)$ is a convex function from $]-\infty,\infty[$ to
$[0,\infty]$ with $\phi(1)=0$. We will consider only
$\phi$-divergences for which the function $\phi(\cdot)$ is strictly
convex and satisfies: the domain of $\phi(\cdot)$, ${\rm dom}\phi:=\{x
\in \mathbb{R}: \phi(x)<\infty\}$ is an interval with end points
$$a_{\phi}<1<b_{\phi},~~ \phi(a_\phi)=\lim_{x\downarrow
a_\phi}\phi(x)~~ \mbox{ and }~~ \phi(a_\phi)=\lim_{x\uparrow b_\phi}\phi(x).$$
The Kullback-Leibler, modified Kullback-Leibler, $\chi^2$,
modified $\chi^2$ and Hellinger  divergences are examples of
$\phi$-divergences; they are obtained respectively for
$\phi(x)=x\log x-x+1$, $\phi(x)=-\log x+x-1$,
$\phi(x)=\frac{1}{2}(x-1)^2$,
$\phi(x)=\frac{1}{2}\frac{(x-1)^2}{x}$ and
$\phi(x)=2(\sqrt{x}-1)^2$. 
The squared Le Cam distance (sometimes called the Vincze-Le Cam distance) and $\mathcal{L}_{1}$-error are obtained respectively for $$\phi(x)=(x-1)^{2}/(2(x-1)) ~~\mbox{and}~~ \phi(x)=|x-1|.$$
We extend the definition of these
divergences on the whole space of all bounded signed  measures via
the extension of the definition of the corresponding $\phi(\cdot)$
functions on the whole real space $\mathbb{R}$ as follows: when
$\phi(\cdot)$ is not well defined on $\mathbb{R}_-$ or well defined but
not convex on $\mathbb{R}$, we set $\phi(x)=+\infty$ for all
$x<0$. Notice that for the $\chi^2$-divergence, the corresponding
$\phi(\cdot)$ function is defined on whole $\mathbb{R}$ and strictly
convex. 
All the above examples are particular cases of the so-called
``\emph{power divergences}'', introduced by
 \cite{CressieRead1984} (see
also \cite[Chapter 2]{LieseVajda1987} and also the \cite{Reny1961}'s paper is to be mentioned here), which are defined through the class of convex real valued
functions, for $\gamma$ in $\mathbb{R}\backslash\left\{0,1\right\}$,
\begin{equation} \label{powerdivergence}
x\in\mathbb{R}_{+}^*\rightarrow
\phi_{\gamma}(x):=\frac{x^{\gamma }-\gamma x+\gamma -1}{\gamma
(\gamma -1)}, 
\end{equation} 
$\phi_{0}(x):=-\log x+x-1$ and $\phi_{1}(x):=x\log x-x+1$. (For
all $\gamma\in\mathbb{R}$, we define
$\phi_\gamma(0):=\lim_{x\downarrow 0}\phi_\gamma (x)$). So, the
$KL$-divergence is associated to $\phi_1$, the $KL_m$ to $\phi_0$,
the $\chi^2$ to $\phi_2$, the $\chi^2_m$ to $\phi_{-1}$ and the
{H}ellinger distance to $\phi_{1/2}$. In the monograph by \cite{LieseVajda1987}
the reader may find detailed ingredients of the
modeling theory as well as surveys of the commonly used divergences.

\noindent Let $\{\mathbb{P}_{\boldsymbol{\theta}}: \boldsymbol{\theta}\in \mathbf{\Theta}\}$ be some identifiable
parametric model with $\mathbf{\Theta}$ a compact subset of $\mathbb{R}^{d}$.
Consider the problem of estimation of the unknown true value of
the parameter $\boldsymbol{\theta}_{0}$ on the basis of an i.i.d. sample
$\mathbf{ X}_{1}, \dots, \mathbf{ X}_{n}$. We shall 
assume that the observed data are from the probability space
$(\mathcal{X}, \mathcal{A}, \mathbb{P}_{\boldsymbol{\theta}_0})$.
\noindent Let $\phi(\cdot)$ be a function of class $\mathcal{C}^2$, strictly convex such that
\begin{equation}
\int \left| \phi ^{\prime }\left( \frac{\mathrm{d}\mathbb{P}_{\boldsymbol{\theta} }(\mathbf{ x})}{\mathrm{d}\mathbb{P}_{\boldsymbol{\alpha}
}(\mathbf{ x})}\right) \right| ~\mathrm{d}\mathbb{P}_{\boldsymbol{\theta} }(\mathbf{ x})<\infty , \forall \boldsymbol{\alpha}\in\mathbf{\Theta}.  \label{condition
integrabilite}
\end{equation}As it is mentioned in \cite{BroniatowskiKeziou2009}, if the function $\phi(\cdot)$ satisfies the following conditions 
\begin{equation}\label{tresPOA}
\begin{array}{l}
 \mbox{ there exists } ~0<\delta<1~ \mbox{ such that for all }~c~ \mbox{ in } ~~[1-\delta,1+\delta],\\   
 \mbox{ we can find numbers }~c_{1},c_{2},c_{3}~\mbox{such that}\\
 \phi(cx)\leq c_{1}\phi(x)+c_{2}|x|+c_{3},~ \mbox{ for all real }~x,       
\end{array}
\end{equation}
then the assumption (\ref{condition
integrabilite}) is satisfied whenever $D_{\phi}(\boldsymbol{\theta},\boldsymbol{\alpha})<\infty$, where $D_{\phi}(\boldsymbol{\theta},\boldsymbol{\alpha})$ stands for the $\phi$-divergence between  $\mathbb{P}_{\boldsymbol{\theta}}$ and $\mathbb{P}_{\boldsymbol{\alpha}}$, refer to  \cite[Lemma 3.2]{BroniatowskiKeziou2006}. Also the real convex functions $\phi(\cdot)$ (\ref{powerdivergence}), associated with the class of power divergences, all satisfy the condition (\ref{condition
integrabilite}), including all standard divergences.  Under assumption (\ref{condition integrabilite}), using Fenchel duality technique,
the divergence $D_{\phi}(\boldsymbol{\theta},\boldsymbol{\theta}_0)$ can be represented
as resulting from an optimization procedure, this result was elegantly proved in \cite{Keziou2003}, \cite{LieseVajda2006} and 
 \cite{BroniatowskiKeziou2009}. \cite{BroniatowskiKeziou2006} called it the
  \emph{dual} form of a divergence, due to its connection with convex analysis.
 According to \cite{LieseVajda2006},  under the strict convexity and the differentiability of the
function $\phi(\cdot)$, it holds

\begin{equation}\label{Convexity}
\phi(t)\geq\phi(s)+\phi^{\prime}(s)(t-s),
\end{equation}
where the equality holds only for $s=t$. Let $\boldsymbol{\theta}$ and $\boldsymbol{\theta}_0$ be fixed
and put $t=\mathrm{d}\mathbb{P}_{\boldsymbol{\theta}}(\mathbf{ x})/\mathrm{d}\mathbb{P}_{\boldsymbol{\theta}_0}(\mathbf{ x})$ and
$s=\mathrm{d}\mathbb{P}_{\boldsymbol{\theta}}(\mathbf{ x})/\mathrm{d}\mathbb{P}_{\boldsymbol{\alpha}}(\mathbf{ x})$ in (\ref{Convexity}) and then
integrate with respect to $\mathbb{P}_{\boldsymbol{\theta}_0}$, to obtain 
\begin{equation}\label{Dualrepresentation}
D_{\phi}(\boldsymbol{\theta}, \boldsymbol{\theta}_{0}):=\int
\phi\left(\frac{\mathrm{d}\mathbb{P}_{\boldsymbol{\theta}}}{\mathrm{d}\mathbb{P}_{\boldsymbol{\theta}_0}}\right)
~\mathrm{d}\mathbb{P}_{\boldsymbol{\theta}_0}=\sup_{\boldsymbol{\alpha}\in
\mathbf{\Theta}}\int h(\boldsymbol{\theta},\boldsymbol{\alpha})~\mathrm{d}\mathbb{P}_{\boldsymbol{\theta}_{0}},
\end{equation}
where $h(\boldsymbol{\theta},\boldsymbol{\alpha},\cdot):\mathbf{ x}\mapsto h(\boldsymbol{\theta},\boldsymbol{\alpha},\mathbf{ x})$ and
\begin{equation}\label{Definition-h}
h(\boldsymbol{\theta},\boldsymbol{\alpha},\mathbf{ x}):=\int \phi ^{\prime }\left( \frac{\mathrm{d}\mathbb{P}_{\boldsymbol{\theta} }}{\mathrm{d}\mathbb{P}_{\boldsymbol{\alpha}
}}\right) ~\mathrm{d}\mathbb{P}_{\boldsymbol{\theta} }-\left[ \frac{\mathrm{d}\mathbb{P}_{\boldsymbol{\theta} }(\mathbf{ x})}{\mathrm{d}\mathbb{P}_{\boldsymbol{\alpha}
}(\mathbf{ x})}\phi ^{\prime }\left( \frac{\mathrm{d}\mathbb{P}_{\boldsymbol{\theta} }(\mathbf{ x})}{\mathrm{d}\mathbb{P}_{\boldsymbol{\alpha} }(\mathbf{ x})}\right)
-\phi \left( \frac{\mathrm{d}\mathbb{P}_{\boldsymbol{\theta} } (\mathbf{ x})}{\mathrm{d}\mathbb{P}_{\boldsymbol{\alpha} }(\mathbf{ x})}\right) \right].
\end{equation}
Furthermore, the supremum in this display (\ref{Dualrepresentation}) is unique and reached in
$\boldsymbol{\alpha}=\boldsymbol{\theta}_0$, independently upon the value of $\boldsymbol{\theta}$.
Naturally, a class of estimators of $\boldsymbol{\theta}_{0}$, called ``dual
$\phi$-divergence estimators'' (D$\phi$DE's), is defined by
\begin{equation}\label{dualestimator}
\widehat{\boldsymbol{\alpha}}_{\phi}(\boldsymbol{\theta}):=\arg\sup_{\boldsymbol{\alpha}\in \mathbf{\Theta}}\mathbb{P}_n
h(\boldsymbol{\theta},\boldsymbol{\alpha}),\;\;\boldsymbol{\theta}\in \mathbf{\Theta},
\end{equation}
where $h(\boldsymbol{\theta},\boldsymbol{\alpha})$ is the function defined in (\ref{Definition-h}) and,
 for a measurable function $f(\cdot)$, $$\mathbb{P}_{n}f:=n^{-1}\sum_{i=1}^{n}f(\mathbf{ X}_i).$$
The class of estimators $\widehat{\boldsymbol{\alpha}}_{\phi}(\boldsymbol{\theta})$ satisfies
\begin{equation}\label{soldualestimator}
\mathbb{P}_n
\frac{\partial}{\partial\boldsymbol{\alpha}}h(\boldsymbol{\theta},\widehat{\boldsymbol{\alpha}}_{\phi}(\boldsymbol{\theta}))=0.
\end{equation}

\noindent Formula (\ref{dualestimator}) defines a family of $M$-estimators indexed by the function $\phi(\cdot)$
specifying the divergence and by some instrumental value of the parameter $\boldsymbol{\boldsymbol{\theta}}$. The $\phi$-divergence estimators are motivated by the fact that a
suitable choice of the divergence may lead to an estimate more
robust than the maximum likelihood estimator (MLE) one, see \cite{JimenezShao2001}.  \cite{TomaBroniatowski2010} studied the robustness of the D$\phi$DE's through the influence function approach, they treated numerous examples of location-scale models and give sufficient conditions for the robustness of D$\phi$DE's. We recall that the maximum likelihood estimate belongs to the
class of estimates (\ref{dualestimator}). Indeed, it is obtained when
$\phi (x)=-\log x+x-1$, that is as the dual modified
$ KL_{m}$-divergence estimate. Observe that 
$\phi^\prime(x)=-\frac{1}{x}+1$
and 
$x  \phi^\prime(x)-\phi(x)=\log x,$ hence $$\int h(\boldsymbol{\theta} ,\boldsymbol{\alpha} ){\rm{d}}\mathbb{P}_{n}=-\int \log
\left( \frac{{\rm{d}}\mathbb{P}_{\boldsymbol{\theta} }}{{\rm{d}}\mathbb{P}_{\boldsymbol{\alpha} }}\right){\rm{d}}\mathbb{P}_{n}.$$
Keeping in mind definitions (\ref{dualestimator}), we get
\begin{eqnarray*}
\widehat{\boldsymbol{\alpha} }_{\rm KL_{m}}(\boldsymbol{\theta} )&=&\arg \sup_{\boldsymbol{\alpha}}-\int \log
\left( \frac{{\rm{d}}\mathbb{P}_{\boldsymbol{\theta} }}{{\rm{d}}\mathbb{P}_{\boldsymbol{\alpha} }}\right) {\rm{d}}\mathbb{P}_{n}\\&=&\arg \sup_{\boldsymbol{\alpha} }\int \log ({\rm{d}}\mathbb{P}_{\boldsymbol{\alpha} }){\rm{d}}\mathbb{P}_{n}={\rm{MLE}},
\end{eqnarray*}
independently upon $\boldsymbol{\theta}$. 

\section{Asymptotic properties}\label{sec2}
In this section, we shall establish the consistency of bootstrapping
under general conditions in the framework of
dual divergence estimation. Define, for a measurable function $f(\cdot)$,
$$\mathbb{P}_n^{\ast}f:=\frac{1}{n}\sum_{i=1}^{n}W_{ni}f(\mathbf{ X}_i),$$ where
$W_{ni}$'s are the bootstrap weights defined on the probability
space $(\mathcal{W},\Omega,\mathbb{P}_{W})$. In view of (\ref{dualestimator}), the
bootstrap estimator can be rewritten as
\begin{equation}\label{Bdualestimator}
\widehat{\boldsymbol{\alpha}}^{\ast}_{\phi}(\boldsymbol{\theta})
:=\arg\sup_{\boldsymbol{\alpha}\in\mathbf{\Theta}}\mathbb{P}_n^{\ast}
h(\boldsymbol{\theta},\boldsymbol{\alpha}).
\end{equation}
The definition of $\widehat{\boldsymbol{\alpha}}^{\ast}_{\phi}(\boldsymbol{\theta})$, defined in (\ref{Bdualestimator}), implies that
\begin{equation}\label{solBdualestimator}
\mathbb{P}_n^{\ast}
\frac{\partial}{\partial\boldsymbol{\alpha}}h(\boldsymbol{\theta},\widehat{\boldsymbol{\alpha}}^{\ast}_{\phi}(\boldsymbol{\theta}))=0.
\end{equation}
The bootstrap weights $W_{ni}$'s are assumed to belong to the
class of exchangeable bootstrap weights introduced in \cite{PraestgaardWellner1993}. In the sequel,  the transpose of a vector $\mathbf{x}$  will be denoted by $\mathbf{x}^\top$.
We shall assume the following conditions.
\begin{enumerate}
\item[W.1] {\it The vector $W_{n}=(W_{n1},\ldots,W_{nn})^{\top}$ is
exchangeable for all $n=1,2,\ldots$, i.e., for any permutation
$\pi=(\pi_{1},\ldots,\pi_{n})$ of $(1,\ldots,n)$, the joint
distribution of $\pi(W_{n})=(W_{n\pi_{1}},\ldots,W_{n\pi_{n}})^{\top}$
is the same as that of $W_{n}$.}

\item[W.2] {\it $W_{ni}\geq 0$ for all $n$, $i$ and
$\sum_{i=1}^{n}W_{ni}=n  $ for all $n$.}

\item[W.3] {\it $\lim\sup_{n\rightarrow\infty}\|W_{n1}\|_{2,1}\leq
C<\infty$, where}
$$\|W_{n1}\|_{2,1}=\int_{0}^{\infty}\sqrt{\mathbb{P}_W(W_{n1}\geq u)}du.$$

\item[W.4] 
$$\lim_{\lambda\rightarrow\infty}\lim\sup_{n\rightarrow\infty}
\sup_{t\geq\lambda}t^{2}\mathbb{P}_W(W_{n1}>t)=0.$$

\item[W.5] {\it
$(1/n)\sum_{i=1}^{n}(W_{ni}-1)^{2}\overset{\mathbb{P}_W}{\longrightarrow}c^{2}>0$.
}
\end{enumerate}
In Efron's nonparametric bootstrap, the bootstrap sample is drawn
from the nonparametric estimate of the true distribution, i.e.,
empirical distribution. Thus, it is easy to show that $W_n\sim
\mbox{Multinomial} (n;n^{-1},\ldots,n^{-1})$ and conditions
{\rm{W.1--W.5}} are satisfied. In general, conditions
{\rm{W.3-W.5}} are easily satisfied under some moment conditions on
$W_{ni}$, see  \cite[Lemma 3.1]{PraestgaardWellner1993}. In addition to Efron's
nonparametric boostrap, the sampling schemes that satisfy
conditions {\rm{W.1--W.5}}, include \emph{Bayesian bootstrap, Multiplier
bootstrap, Double bootstrap, and Urn boostrap}. This list is sufficiently long to indicate
that conditions {\rm{W.1--W.5}}, are not unduely restrictive. Notice that the value of $c$ in {\rm{W.5}
 is independent of $n$
and depends on the resampling method, e.g., $c = 1$ for the nonparametric
bootstrap and Bayesian bootstrap, and $c =\sqrt{2}$ for the double bootstrap. A more precise discussion
of this general formulation of the bootstrap can be found in \cite{PraestgaardWellner1993},
\cite{vanderVaartWellner1996} and \cite{Kosorok2008}.

There exist two sources of randomness for the bootstrapped quantity, i.e.,
$\widehat{\boldsymbol{\alpha}}^{\ast}_{\phi}(\boldsymbol{\theta})$: the first comes from
the observed data and the  second is due to the resampling done by the bootstrap, i.e.,
random $W_{ni}$'s. Therefore, in order to rigorously state our
main theoretical results for the general bootstrap of $\phi$-divergence estimates, we need to
specify relevant probability spaces and define stochastic orders
with respect to relevant probability measures. Following \cite{ChengHuang2010} and \cite{WellnerZhan1996},  
we shall view $\mathbf{ X}_i$ as the $i$-th coordinate projection from the canonical probability space $(\mathcal{X}^{\infty},
\mathcal{A}^{\infty}, \mathbb{P}_{\boldsymbol{\theta}_0}^{\infty})$ onto the $i$-th copy of $\mathcal{X}$.
For the joint randomness involved, the product probability space is defined as
$$
(\mathcal{X}^{\infty}, \mathcal{A}^{\infty},\mathbb{P}_{\boldsymbol{\theta}_0}^{\infty})\times(\mathcal{W},\Omega, \mathbb{P}_W)=
(\mathcal{X}^{\infty}\times\mathcal{W},\mathcal{A}^{\infty}\times \Omega, \mathbb{P}_{\boldsymbol{\theta}_0}^{\infty}\times \mathbb{P}_W).
$$
Throughout the
paper, we assume that the bootstrap weights $W_{ni}$'s are
independent of the data $\mathbf{ X}_i$'s, thus $$\mathbb{P}_{XW}=\mathbb{P}_{\boldsymbol{\theta}_0}\times \mathbb{P}_W.$$
Given a real-valued function $\Delta_n$ defined on the above product
probability space, e.g. $\widehat{\boldsymbol{\alpha}}^{\ast}_{\phi}(\boldsymbol{\theta})$, we say that
$\Delta_n$ is of an order $o^{o}_{\mathbb{P}_W}(1)$ in
$\mathbb{P}_{\boldsymbol{\theta}_0}$-probability if, for any $\epsilon,\eta>0$, as $n\rightarrow 0$,
\begin{eqnarray}
\mathbb{P}_{\boldsymbol{\theta}_0}\{P^{o}_{W|X}(|\Delta_n|>\epsilon)>\eta\}&\longrightarrow&
0,
\end{eqnarray}
 and that $\Delta_n$ is of an order
$O^{o}_{\mathbb{P}_W}(1)$ in $\mathbb{P}_{\boldsymbol{\theta}_0}$-probability if, for any $\eta>0$,
there exists a $0<M<\infty$ such that, as $n\rightarrow 0$,
\begin{eqnarray}
\mathbb{P}_{\boldsymbol{\theta}_0}\{P^{o}_{W|X}(|\Delta_n|\geq M)> \eta\}\longrightarrow
0,
\end{eqnarray} where the superscript ``$o$'' denotes the outer probability, see \cite{vanderVaartWellner1996} for more details on outer probability measures. For more details on stochastic orders, the interested reader may refer to \cite{ChengHuang2010}, in particular, Lemma 3 of the cited reference.\\
To establish the consistency of $\widehat{\boldsymbol{\alpha}}^{\ast}_{\phi}(\boldsymbol{\theta})$,
the following conditions are assumed in our analysis.
\begin{enumerate}
\item[(A.1)]
\begin{equation}\label{well-sep}
  \mathbb{P}_{\boldsymbol{\theta}_0}h(\boldsymbol{\theta},\boldsymbol{\theta}_{0})>\sup_{\boldsymbol{\alpha}\not \in N(\boldsymbol{\theta}_{0})}\mathbb{P}_{\boldsymbol{\theta}_0}
h(\boldsymbol{\theta},\boldsymbol{\alpha})
\end{equation}
for any open set $N(\boldsymbol{\theta}_{0})\subset\mathbf{\Theta}$ containing
$\boldsymbol{\theta}_{0}$.
\item[(A.2)]\begin{equation}
\sup_{\boldsymbol{\alpha}\in\mathbf{\Theta}}|\mathbb{P}_n^{\ast}
h(\boldsymbol{\theta},\boldsymbol{\alpha})-\mathbb{P}_{\boldsymbol{\theta}_0} h(\boldsymbol{\theta},\boldsymbol{\alpha})|\overset{\mathbb{P}_{XW}^{o}}
{\longrightarrow}0.\label{concon}
\end{equation}
\end{enumerate}
The following theorem gives the consistency of the bootstrapped
estimates
$\widehat{\boldsymbol{\alpha}}^{\ast}_{\phi}(\boldsymbol{\theta})$.
\begin{theorem}\label{asythm-a} Assume that conditions {\rm{(A.1)}} and {\rm{(A.2)}}
hold.  Suppose that
conditions {\rm{(A.3--5)}} and {\rm{W.1--W.5}} hold.
Then
$\widehat{\boldsymbol{\alpha}}^{\ast}_{\phi}(\boldsymbol{\theta})$
is a consistent estimate of $\boldsymbol{\theta}_0$. That is
$$\widehat{\boldsymbol{\alpha}}^{\ast}_{\phi}(\boldsymbol{\theta})\overset{\mathbb{P}_{W}^{o}}{\longrightarrow}\boldsymbol{\theta}_0
~~\mbox{ in }~~ \mathbb{P}_{\boldsymbol{\theta}_0}\mbox{-probability}.$$
\end{theorem}
\noindent The proof of Theorem \ref{asythm-a} is postponed  until \S \ref{proof}.\\ We need the following definitions, refer to \cite{vanderVaart1998} and \cite{vanderVaartWellner1996} among others.
If $\mathcal{F}$ is a class of functions for which, we have almost surely,
$$
\|\mathbb{P}_n- \mathbb{P}\|_\mathcal{F} = \sup_{f\in\mathcal{F}}
|\mathbb{P}_nf-\mathbb{P}f| \rightarrow 0,
$$
then we say that $\mathcal{F}$ is a $\mathbb{P}$-Glivenko-Cantelli class of functions.
 If $\mathcal{F}$ is a class of functions
for which
$$
\mathbb{G}_n=\sqrt{n}(\mathbb{P}_n-\mathbb{P})\rightarrow \mathbb{G} ~~\mbox{ in } ~~\ell^\infty(\mathcal{F}),
$$
where $\mathbb{G}$ is a mean-zero $\mathbb{P}$-Brownian bridge process with (uniformly-) continuous sample
paths with respect to the semi-metric $\rho_\mathbb{P} (f, g)$, defined by
$$
\rho_\mathbb{P}^2(f, g) = Var_\mathbb{P} (f(X)-g(X)),
$$
then we say that $\mathcal{F}$ is a $\mathbb{P}$-Donsker class of functions. Here
$$
\ell^\infty(\mathcal{F})=\left\{v: \mathcal{F}\mapsto \mathbb{R}\Big|\|v\|_\mathcal{F}=\sup_{f\in\mathcal{F}}|v(f)|<\infty\right\}
$$
and $\mathbb{G}$ is a $\mathbb{P}$-Brownian bridge process on $\mathcal{F}$ if it is a mean-zero Gaussian process with
covariance function
$$
\mathbb{E}(\mathbb{G}(f)\mathbb{G}(g)) = \mathbb{P}fg - (\mathbb{P}f)(\mathbb{P}g).
$$

\begin{remark}
\begin{itemize}
  \item Condition {\rm{(A.1)}} is the ``well separated" condition, compactness of the parameter
   space $\mathbf{\Theta}$ and the continuity of divergence imply that the optimum is well-separated,
   provided the parametric model is identified, see \cite[Theorem 5.7]{vanderVaart1998}.

  \item Condition {\rm{(A.2)}} holds if the class $$\left\{h(\boldsymbol{\theta},\boldsymbol{\alpha}):
  ~\boldsymbol{\alpha}\in\mathbf{\Theta}\right\}$$ is shown to be $\mathbb{P}$-Glivenko-Cantelli, by applying 
    \cite[Lemma 3.6.16]{vanderVaartWellner1996} and   \cite[Lemma A.1]{ChengHuang2010}.
\end{itemize}
\end{remark}

\noindent For any fixed $\delta_n>0$, define the class of
functions $\mathcal{H}_{n}$ and $\dot{\mathcal{H}}_{n}$ as
\begin{equation}\label{H}
\mathcal{H}_{n}:=\left\{\frac{\partial}{\partial\boldsymbol{\alpha}}h(\boldsymbol{\theta},
\boldsymbol{\alpha}):~\|\boldsymbol{\alpha}-\boldsymbol{\theta}_0\|\leq\delta_n\right\}
\end{equation}
and
\begin{equation}\label{dotH}
\dot{\mathcal{H}}_{n}:=\left\{\frac{\partial^2}{\partial\boldsymbol{\alpha}^2}
h(\boldsymbol{\theta},\boldsymbol{\alpha}):~\|\boldsymbol{\alpha}-\boldsymbol{\theta}_0\|\leq\delta_n\right\}.
\end{equation}
We shall say a class of functions
$\mathcal{H}\in M(\mathbb{P}_{\boldsymbol{\theta}_0})$ if $\mathcal{H}$ possesses enough
measurability for randomization with i.i.d. multipliers to be possible, i.e., $\mathbb{P}_{n}$
can be randomized, in other word, we
can replace $(\delta_{\mathbf{ X}_{i}}-\mathbb{P}_{\boldsymbol{\theta}_0})$ by $(W_{ni}-1)\delta_{\mathbf{ X}_{i}}$.
It is known that $\mathcal{H}\in M(\mathbb{P}_{\boldsymbol{\theta}_0})$, e.g., if $\mathcal{H}$
 is countable, or if $\{\mathbb{P}_{n}\}_{n}^{\infty}$ are stochastically separable in  $\mathcal{H}$, or if $\mathcal{H}$ is image admissible Suslin; see \cite[pages 853 and 854]{GineZinn1990}.\\
To state our result concerning the asymptotic normality, we shall assume the following additional conditions.
\begin{enumerate}
 \item[(A.3)] The matrices $$\displaystyle{V:=\mathbb{P}_{\boldsymbol{\theta}_0}\frac{\partial}{\partial\boldsymbol{\alpha}}h(\boldsymbol{\theta},\boldsymbol{\theta}_0)\frac{\partial}{\partial\boldsymbol{\alpha}}h(\boldsymbol{\theta},\boldsymbol{\theta}_0)^{\top}}$$
and $$\displaystyle{S:=-\mathbb{P}_{\boldsymbol{\theta}_0}\frac{\partial^2}{\partial\boldsymbol{\alpha}^2}h(\boldsymbol{\theta},\boldsymbol{\theta}_0)}$$ are
non singular.
 \item[(A.4)] The class $\displaystyle{\mathcal{H}_{n}\in M(\mathbb{P}_{\boldsymbol{\theta}_0})\cap
L_{2}(\mathbb{P}_{\boldsymbol{\theta}_0})}$ and is $\mathbb{P}$-Donsker.
\item[(A.5)] The class $\displaystyle{\dot{\mathcal{H}}_{n}\in M(\mathbb{P}_{\boldsymbol{\theta}_0})\cap
L_{2}(\mathbb{P}_{\boldsymbol{\theta}_0})}$ and is $\mathbb{P}$-Donsker.
\end{enumerate}
Conditions {\rm{(A.4)}} and {\rm{(A.5)}}  ensure that the ``size''
of the function classes $\mathcal{H}_{n}$ and $\dot{\mathcal{H}}_{n}$ are reasonable so that
the bootstrapped empirical processes
$$\mathbb{G}_n^{\ast}\equiv\sqrt{n}(\mathbb{P}_n^{\ast}-\mathbb{P}_n)$$
indexed, respectively by $\mathcal{H}_{n}$ and $\dot{\mathcal{H}}_n$, have a limiting process conditional
on the original observations, we refer for instance to \cite[Theorem 2.2]{PraestgaardWellner1993}.
The main result to be proved here may now be stated precisely as
follows.
\begin{theorem}\label{asythm-b}
Assume that
$\widehat{\boldsymbol{\alpha}}_{\phi}(\boldsymbol{\theta})$ and
$\widehat{\boldsymbol{\alpha}}^{\ast}_{\phi}(\boldsymbol{\theta})$
fullfil (\ref{soldualestimator}) and (\ref{solBdualestimator}),
respectively. In addition suppose that
$$\widehat{\boldsymbol{\alpha}}_{\phi}(\boldsymbol{\theta})\overset
{\mathbb{P}_{\boldsymbol{\theta}_0}}{\longrightarrow}\boldsymbol{\theta}_0
~~\mbox{ and }~~
\widehat{\boldsymbol{\alpha}}^{\ast}_{\phi}(\boldsymbol{\theta})
\overset{\mathbb{P}_{W}^{o}}{\longrightarrow}\boldsymbol{\theta}_0~~
\mbox{in}~~ \mathbb{P}_{\boldsymbol{\theta}_0}\mbox{-probability.}$$ Assume that
conditions {\rm{(A.3--5)}} and {\rm{W.1--W.5}} hold.
 Then we have
\begin{eqnarray}
\|\widehat{\boldsymbol{\alpha}}^{\ast}_{\phi}(\boldsymbol{\theta})-\boldsymbol{\theta}_{0}\|=O^{o}_{\mathbb{P}_{W}}(n^{-1/2})
\label{bconratep}
\end{eqnarray}
in $\mathbb{P}_{\boldsymbol{\theta}_0}$-probability. Furthermore,
\begin{eqnarray}\label{bcons}
\sqrt{n}(\widehat{\boldsymbol{\alpha}}^{\ast}_{\phi}(\boldsymbol{\theta})-\widehat{\boldsymbol{\alpha}}_{\phi}(\boldsymbol{\theta}))=
-S^{-1}\mathbb{G}_{n}^{\ast}\frac{\partial}{\partial\boldsymbol{\alpha}}h(\boldsymbol{\theta},\boldsymbol{\theta}_0)+
o^{o}_{\mathbb{P}_{W}}(1)
\end{eqnarray}
in $\mathbb{P}_{\boldsymbol{\theta}_0}$-probability. Consequently,
\begin{eqnarray}\label{bconcor}
\sup_{\mathbf{ x}\in\mathbb{R}^d}\left|\mathbb{P}_{W|\mathcal{X}_n}((\sqrt{n}/c)(\widehat{\boldsymbol{\alpha}}^{\ast}_{\phi}
(\boldsymbol{\theta})-\widehat{\boldsymbol{\alpha}}_{\phi}(\boldsymbol{\theta}))\leq \mathbf{ x})-\mathbb{P}(N(0,\Sigma)\leq
\mathbf{ x})\right|=o_{\mathbb{P}_{\boldsymbol{\theta}_0}}(1),
\end{eqnarray}
where ``$\leq$'' is taken
componentwise and ``$c$'' is given in {\rm{W.5}}, whose value depends on the used sampling
scheme, and $$\Sigma\equiv S^{-1}V (S^{-1})^{\top}$$ where $S$ and $V$
are given in condition {\rm{(A.3)}}. Thus, we have
\begin{eqnarray}
\sup_{\mathbf{ x}\in\mathbb{R}^d}\left|\mathbb{P}_{W|\mathcal{X}_n}((\sqrt{n}/c)(\widehat{\boldsymbol{\alpha}}^{\ast}_{\phi}(\boldsymbol{\theta})-\widehat{\boldsymbol{\alpha}}_{\phi}(\boldsymbol{\theta}))\leq \mathbf{ x})-\mathbb{P}_{\boldsymbol{\theta}_0}(\sqrt{n}(\widehat{\boldsymbol{\alpha}}_{\phi}(\boldsymbol{\theta})-\boldsymbol{\theta}_0)\leq
\mathbf{ x})\right|\overset{\mathbb{P}_{\boldsymbol{\theta}_0}}{\longrightarrow}0\label{proconv}
\end{eqnarray}
\end{theorem}
\noindent The proof of Theorem \ref{asythm-a} is captured in the forthcoming \S \ref{proof}.\\
\begin{remark}
Note that an appropriate choice of the the bootstrap weights $W_{ni}$'s implicates a smaller limit variance, that is, $c^{2}$
is smaller than $1$. For instance, typical examples are i.i.d.-weighted bootstraps and the multivariate hypergeometric bootstrap, refer to    
\cite[Examples 3.1 and 3.4]{PraestgaardWellner1993}.
\end{remark}
\noindent Following \cite{ChengHuang2010}, we shall illustrate how to apply our results to construct the confidence sets.
A lower $\epsilon$-th quantile of bootstrap distribution is defined to be any  $ q_{n\epsilon}^{\ast}\in\mathbb{R}^d$
fulfilling
$$ q_{n\epsilon}^{\ast}:=\inf\{\mathbf{ x}:
\mathbb{P}_{W|\mathcal{X}_n}(\widehat{\boldsymbol{\alpha}}^{\ast}_{\phi}(\boldsymbol{\theta})\leq \mathbf{ x})\geq\epsilon\},$$
where $\mathbf{ x}$ is an infimum over the given set only if there
does not exist a $\mathbf{x}_1<\mathbf{ x}$ in $\mathbb{R}^d$ such that
$$\mathbb{P}_{W|\mathcal{X}_n}(\widehat{\boldsymbol{\alpha}}^{\ast}_{\phi}(\boldsymbol{\theta})\leq \mathbf{x}_1)\geq\epsilon.$$
Keep in mind the assumed regularity conditions on the criterion function, that is,
$h(\boldsymbol{\theta},\boldsymbol{\alpha})$ in the  present framework, we can, without loss of
generality, suppose that $$\mathbb{P}_{W|\mathcal{X}_n}(\widehat{\boldsymbol{\alpha}}^{\ast}_{\phi}(\boldsymbol{\theta})
\leq q_{n\epsilon}^{\ast})=\epsilon.$$ Making use the distribution consistency result given in
(\ref{proconv}), we can approximate the $\epsilon$-th quantile of the distribution of $$(\widehat{\boldsymbol{\alpha}}_{\phi}(\boldsymbol{\theta})-\boldsymbol{\theta}_0)~~\mbox{ by }~~(q_{n\epsilon}^{\ast}-\widehat{\boldsymbol{\alpha}}_{\phi}(\boldsymbol{\theta}))/c.$$ Therefore, we define the \emph{percentile}-type bootstrap confidence set as 
\begin{equation}\label{Cpercentile}
\mathrm{C}(\epsilon):= \left[\widehat{\boldsymbol{\alpha}}_{\phi}(\boldsymbol{\theta})+\frac{ q_{n(\epsilon/2)
}^{\ast}-\widehat{\boldsymbol{\alpha}}_{\phi}(\boldsymbol{\theta})} {c},\widehat{\boldsymbol{\alpha}}_{\phi}(\boldsymbol{\theta})+\frac{ q_{n(1-\epsilon/2)}
^{\ast}-\widehat{\boldsymbol{\alpha}}_{\phi}(\boldsymbol{\theta})}
{c}\right].
\end{equation}
In a similar manner, the $\epsilon$-th quantile of $\sqrt{n}(\widehat{\boldsymbol{\alpha}}_{\phi}(\boldsymbol{\theta})-\boldsymbol{\theta}_0)$ can be approximated by $\widetilde{q}_{n\epsilon}^{\ast}$, where $\widetilde{q}_{n\epsilon}^{\ast}$ is the $\epsilon$-th quantile of the hybrid quantity 
$(\sqrt{n}/c)(\widehat{\boldsymbol{\alpha}}^{\ast}_{\phi}(\boldsymbol{\theta})
-\widehat{\boldsymbol{\alpha}}_{\phi}(\boldsymbol{\theta}))$, i.e., 
$$\mathbb{P}_{W|\mathcal{X}_n}((\sqrt{n}/c)(\widehat{\boldsymbol{\alpha}}^{\ast}_{\phi}(\boldsymbol{\theta})
-\widehat{\boldsymbol{\alpha}}_{\phi}(\boldsymbol{\theta}))\leq\widetilde{q}_{n \epsilon}^{\ast})=\epsilon.$$
Note that 
$$\widetilde{q}_{n\epsilon}^{\ast}=(\sqrt{n}/c)( q_{n\epsilon}^{\ast}-
\widehat{\boldsymbol{\alpha}}_{\phi}(\boldsymbol{\theta})).$$
 Thus, the \emph{hybrid}-type bootstrap confidence set would be defined as follows 
\begin{equation}\label{Chybrid}
\widetilde{\mathrm{C}}(\epsilon):=\left[\widehat{\boldsymbol{\alpha}}_{\phi}(\boldsymbol{\theta})-\frac{\widetilde{q}_{n(1-\epsilon/2)}^{\ast}}{\sqrt{n}}
,\widehat{\boldsymbol{\alpha}}_{\phi}(\boldsymbol{\theta})-\frac{\widetilde{q}_{n(\epsilon/2)}^{\ast}}
{\sqrt{n}}\right].
\end{equation}
 Note that
$ q_{n\epsilon}^{\ast}$ and $\widetilde{q}_{n\epsilon}^{\ast}$ are not
unique by the fact that we assume $\boldsymbol{\theta}$ is a vector. Recall that, for any $\mathbf{ x}\in\mathbb{R}^d$,
\begin{eqnarray*}
\mathbb{P}_{\boldsymbol{\theta}_0}(\sqrt{n}(\widehat{\boldsymbol{\alpha}}_{\phi}(\boldsymbol{\theta})-\boldsymbol{\theta}_0)\leq \mathbf{ x})&\longrightarrow&
\Psi(\mathbf{ x}),\\
\mathbb{P}_{W|\mathcal{X}_n}((\sqrt{n}/c)(\widehat{\boldsymbol{\alpha}}^{\ast}_{\phi}(\boldsymbol{\theta})-
\widehat{\boldsymbol{\alpha}}_{\phi}(\boldsymbol{\theta}))\leq \mathbf{ x})&\overset{\mathbb{P}_{\boldsymbol{\theta}_0}}{\longrightarrow}&\Psi(\mathbf{ x}),
\end{eqnarray*}
where $$\Psi(\mathbf{ x})=\mathbb{P}(N(0,\Sigma)\leq \mathbf{ x}).$$ According to  the quantile convergence Theorem,
i.e.,  \cite[Lemma 21.1]{vanderVaart1998}, we have, almost surely, 
$$\widetilde{q}_{n\epsilon}^{\ast}\overset{\mathbb{P}_{XW}}{\longrightarrow}
\Psi^{-1}(\epsilon).$$
When applying quantile convergence theorem, we use the almost sure representation, that is, \cite[Theorem 2.19]{vanderVaart1998}, and argue along subsequences. Considering the Slutsky's Theorem which
ensures that
$$\sqrt{n}(\widehat{\boldsymbol{\alpha}}_{\phi}(\boldsymbol{\theta})-\boldsymbol{\theta}_0)-\widetilde{q}_{n(\epsilon/2)}^{\ast}~~\mbox{
weakly converges to }~~ N(0,\Sigma)-\Psi^{-1}(\epsilon/2),$$ we further have
\begin{eqnarray*}
\mathbb{P}_{XW}\left(\boldsymbol{\theta}_0\leq\widehat{\boldsymbol{\alpha}}_{\phi}(\boldsymbol{\theta})-
\frac{\widetilde{q}_{n(\epsilon/2)}^{\ast}}{\sqrt{n}}\right)&=&\mathbb{P}_{XW}\left(
\sqrt{n}(\widehat{\boldsymbol{\alpha}}_{\phi}(\boldsymbol{\theta})-\boldsymbol{\theta}_0)\geq\widetilde{q}_{n(\epsilon/2)}^{\ast}
\right)\\&&\longrightarrow\mathbb{P}_{XW}\left(N(0,\Sigma)\geq
\Psi^{-1}(\epsilon/2)\right)\\&&=1-\epsilon/2.
\end{eqnarray*}
The above arguments prove the consistency of the {\it hybrid}-type
bootstrap confidence set, i.e., (\ref{perci1}), and can also be
applied to the {\it percentile}-type bootstrap confidence set,
i.e., (\ref{perci2}). For an in-depth study and more rigorous proof, we may refer to \cite[Lemma 23.3]{vanderVaart1998}.
The above discussion may be summarized as follows.
\begin{corollary}\label{perci}
Under the conditions in Theorem~\ref{asythm-b}, we have, as $n\rightarrow\infty$,
\begin{equation}
\mathbb{P}_{XW}\left(\widehat{\boldsymbol{\alpha}}_{\phi}(\boldsymbol{\theta})+\frac{ q_{n(\epsilon/2)
}^{\ast}-\widehat{\boldsymbol{\alpha}}_{\phi}(\boldsymbol{\theta})} {c}\leq\boldsymbol{\theta}_0
\leq\widehat{\boldsymbol{\alpha}}_{\phi}(\boldsymbol{\theta})+\frac{ q_{n(1-\epsilon/2)}
^{\ast}-\widehat{\boldsymbol{\alpha}}_{\phi}(\boldsymbol{\theta})}
{c}\right)\longrightarrow1-\epsilon,\label{perci2}
\end{equation}
\begin{equation}
\mathbb{P}_{XW}\left(\widehat{\boldsymbol{\alpha}}_{\phi}(\boldsymbol{\theta})-\frac{\widetilde{q}_{n(1-\epsilon/2)}^{\ast}}{\sqrt{n}}
\leq\boldsymbol{\theta}_0\leq\widehat{\boldsymbol{\alpha}}_{\phi}(\boldsymbol{\theta})-\frac{\widetilde{q}_{n(\epsilon/2)}^{\ast}}
{\sqrt{n}}\right)\longrightarrow1-\epsilon.\label{perci1}
\end{equation}
\end{corollary}
\noindent It is well known that the above bootstrap confidence sets can be
obtained easily through routine bootstrap sampling.
\begin{remark}
Notice that the choice of weights depends  on the problem at hand : accuracy of the estimation of the entire
 distribution of the  statistic, accuracy of a confidence interval, accuracy in large deviation sense, accuracy for a
  finite sample size, we may refer to \cite{Lancelot1997} and the references therein for more details. \cite{Bertail95} indicate that the area where the weighted  bootstrap clearly performs better
   than the classical bootstrap is in term of coverage accuracy.
\end{remark}
\subsection{On the choice of the escort parameter}\label{choiceoftheescort}

\noindent The very peculiar choice of the escort parameter defined through $\boldsymbol{\theta}=\boldsymbol{\theta}_0$ has same limit properties as the MLE one. The D$\phi$DE $\widehat{\boldsymbol{\alpha}}_\phi\left(\boldsymbol{\theta}_{0}\right)$, in this case, has variance which indeed coincides with the MLE one, see for instance \cite[Theorem 2.2, (1) (b)]{Keziou2003}. This result is of some relevance, since it leaves open the choice of the divergence, while keeping good asymptotic properties.
\noindent For data generated from the distribution $\mathcal{N}(0,1)$, Figure \ref{normallocation1} shows that the global maximum of the empirical criterion $\mathbb{P}_nh\left(\widehat{\boldsymbol{\theta}}_n,\boldsymbol{\alpha}\right)$ is zero, independently of the value of the escort parameter $\widehat{\boldsymbol{\theta}}_n$ (the sample mean $\overline{X}=n^{-1}\sum_{i=1}^{n}\mathbf{ X}_{i}$, in Figure \ref{normallocation1}(a) and the median in Figure \ref{normallocation1}(b)) for all the considered divergences which is in agreement with the result of \cite[Theorem 6]{Broniatowski2011},  
where 
it is showed  that all differentiable divergences produce the same estimator of the parameter on any regular exponential family, in particular the normal models, which is the MLE one, provided that the conditions (\ref{tresPOA}) and $D_{\phi}(\boldsymbol{\theta},\boldsymbol{\alpha})<\infty$ are satisfied.

\begin{figure}[!ht]
\begin{center}
\centerline{\includegraphics[width=8cm]{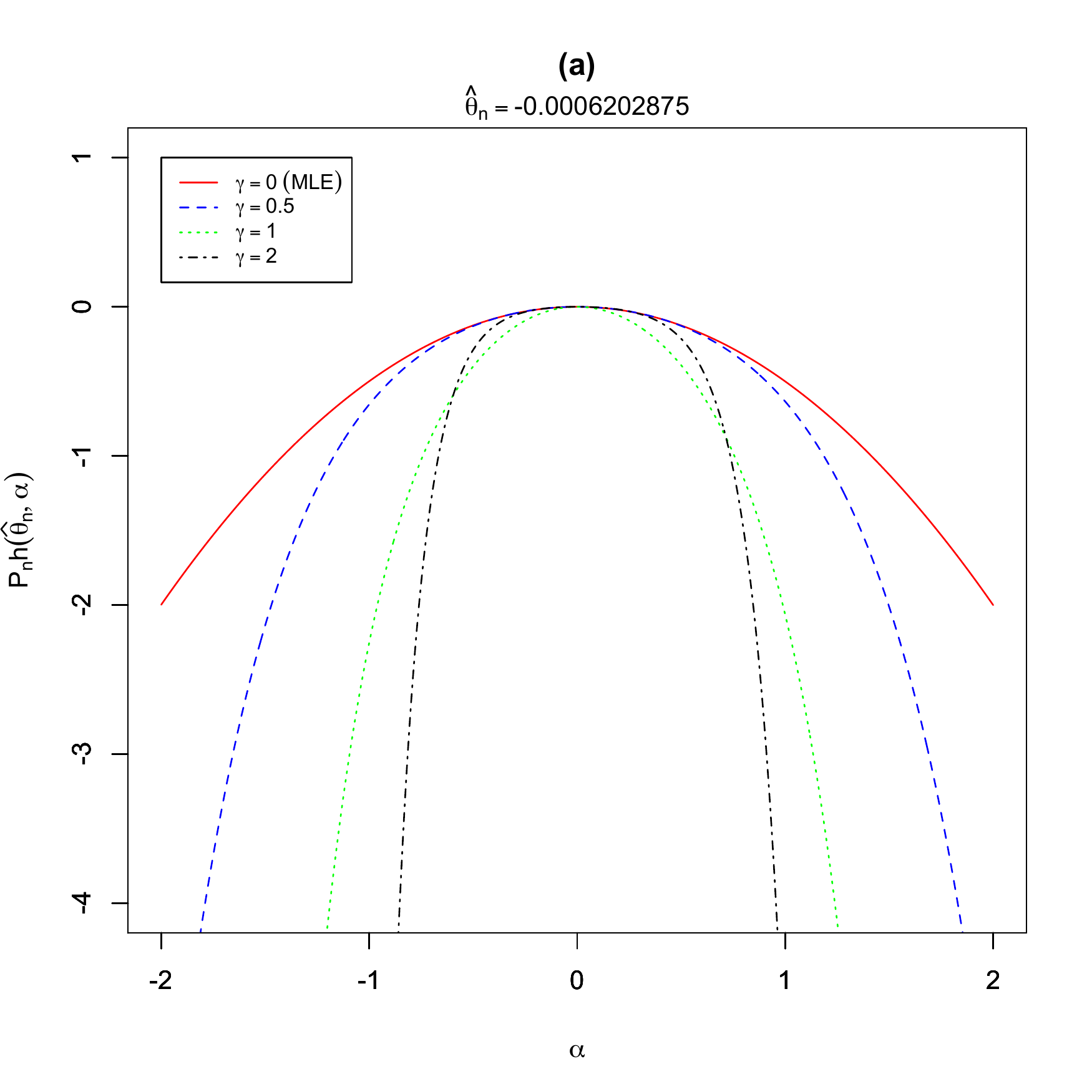}
\includegraphics[width=8cm]{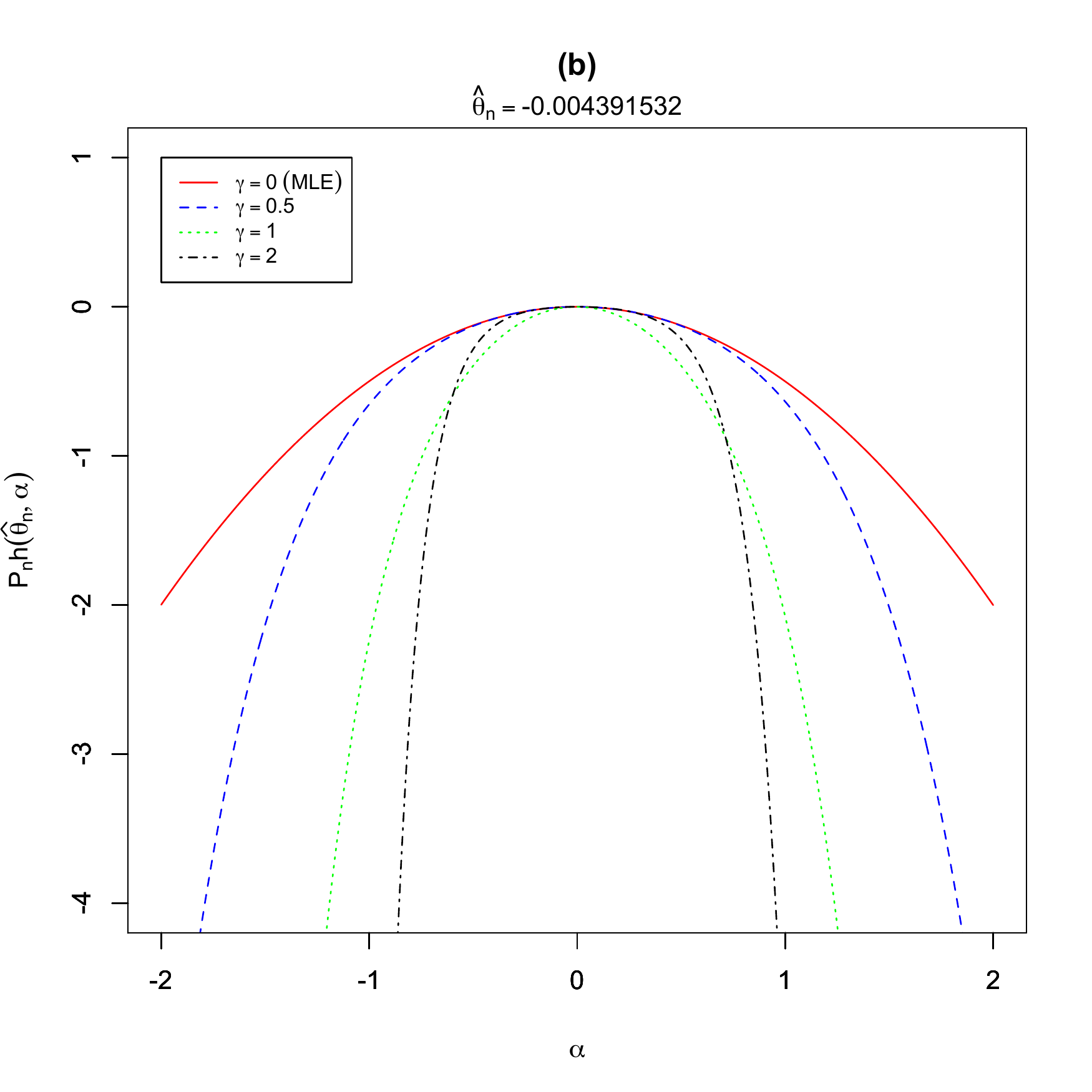}}
\end{center}
\caption{Criterion for the normal location model.} \label{normallocation1}
\end{figure}

\begin{figure}[!ht]
\begin{center}
\centerline{
\includegraphics[width=8cm]{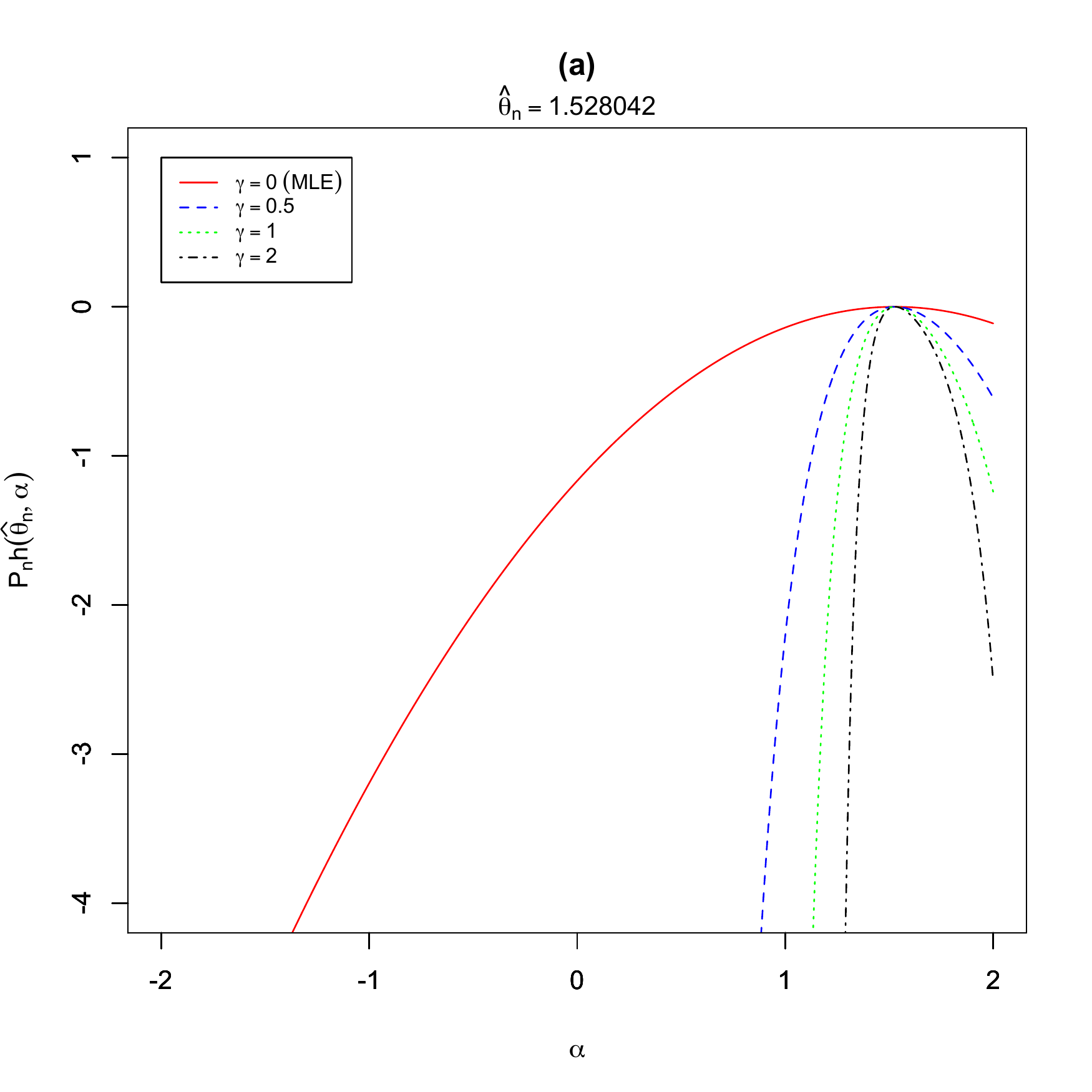}
\includegraphics[width=8cm]{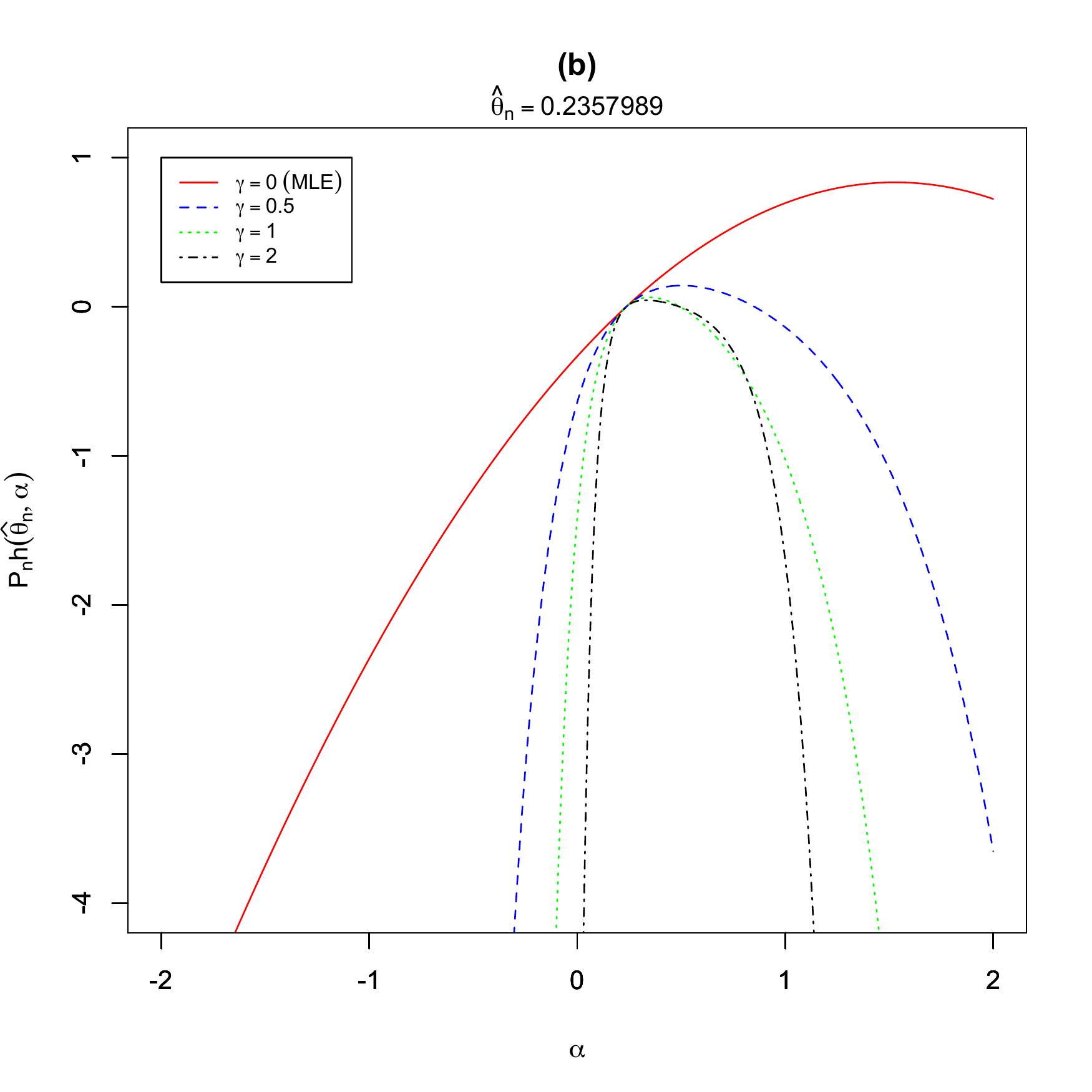}}
\end{center}
\caption{Criterion for the normal location model under contamination.} \label{normallocation2}
\end{figure}

\noindent Unlike the case of data without contamination, the choice of the escort parameter is crucial in the estimation method in the presence of outliers. We plot in Figure \ref{normallocation2} the empirical criterion $\mathbb{P}_nh\left(\widehat{\boldsymbol{\theta}}_n,\boldsymbol{\alpha}\right)$, where the data are generated from the distribution
\begin{equation*}
  (1-\epsilon)\mathcal{N}(\boldsymbol{\theta}_0,1)+\epsilon\delta_{10},
\end{equation*}
where $\epsilon=0.1$, $\boldsymbol{\theta}_0 = 0$ and $\delta_{x}$ stands for the Dirac measure at $x$. Under contamination, when we take the empirical ``\emph{mean}'', $\widehat{\boldsymbol{\theta}}_n=\overline{X}$, as the value of the escort parameter $\boldsymbol{\theta}$, Figure \ref{normallocation2}(a)  shows how the global maximum of the empirical criterion $\mathbb{P}_nh\left(\widehat{\boldsymbol{\theta}}_n,\boldsymbol{\alpha}\right)$ shifts from zero to the contamination point. In Figure \ref{normallocation2}(b), the choice of the ``\emph{median}'' as escort parameter value leads to the position of the global maximum
remains close to $ \boldsymbol{\alpha}= 0$, for {H}ellinger ($\gamma=0.5$), $\chi^2$ ($\gamma=2$)  and $KL$-divergence ($\gamma=1$), while the criterion associated to the $KL_m$-divergence ($\gamma=0$, the maximum is the MLE) stills affected by the presence of outliers.

\noindent In practice, the consequence is that if the data are subject to contamination the escort parameter should be chosen as a robust estimator of $\boldsymbol{\theta}_0$, say $\widehat{\boldsymbol{\theta}}_{n}$. For more details  about the performances of dual
$\phi$-divergence estimators  for normal density models,  we refer to \cite{Cherfi2011a}.

\section{Examples}\label{Examples} Keep in mind the definitions (\ref{Dualrepresentation}) and (\ref{Definition-h}). In what follows, for easy reference and completeness, we give some usual examples of divergences, discussed in \cite{BouzebdaKeziou2010b,BouzebdaKeziou2010a}, of divergences and the
associated estimates, we may refer also to \cite{BroniatowskiVajda2009} for more examples and details.
\begin{itemize}
  \item Our first example  is the  Kullback-Leibler divergence
  \begin{eqnarray*}
  \phi(x)&=&x\log x-x+1\\
\phi^\prime(x)&=&\log x\\
x  \phi^\prime(x)-\phi(x)&=&x-1.
\end{eqnarray*}
The estimate of $D_{\rm KL}(\boldsymbol{\theta},\boldsymbol{\theta}_0)$ is given
by
\begin{eqnarray*}
\widehat{D}_{\rm
KL}(\boldsymbol{\theta},\boldsymbol{\theta}_0)&=&\sup_{\boldsymbol{\alpha}\in \mathbf{\Theta}}\left\{\int\log\left(\frac{{\rm d}\mathbb{P}_{\boldsymbol{\theta}}}{{\rm d}\mathbb{P}_{\boldsymbol{\alpha}}}
\right){\rm d}\mathbb{P}_{\boldsymbol{\theta}}-\int\left(\frac{{\rm d}\mathbb{P}_{\boldsymbol{\theta}}}{{\rm d}\mathbb{P}_{\boldsymbol{\alpha}}}-1
\right){\rm d}\mathbb{P}_n\right\}
\end{eqnarray*}
and the estimate of the parameter $\boldsymbol{\theta}_0$, with escort parameter $\boldsymbol{\theta}$, is defined
as follows
\begin{equation*}
\widehat{\boldsymbol{\alpha}}_{\rm KL}(\boldsymbol{\theta}) := \arg\sup_{\boldsymbol{\alpha}\in \mathbf{\Theta}}
\left\{\int\log\left(\frac{{\rm d}\mathbb{P}_{\boldsymbol{\theta}}}{{\rm d}\mathbb{P}_{\boldsymbol{\alpha}}}
\right){\rm d}\mathbb{P}_{\boldsymbol{\theta}}-\int\left(\frac{{\rm d}\mathbb{P}_{\boldsymbol{\theta}}}{{\rm d}\mathbb{P}_{\boldsymbol{\alpha}}}-1
\right){\rm d}\mathbb{P}_n\right\}.
\end{equation*}
  \item The second one is the $\chi^2$-divergence
  \begin{eqnarray*}
  \phi(x)&=&\frac{1}{2}(x-1)^2\\
\phi^\prime(x)&=&x-1\\
x  \phi^\prime(x)-\phi(x)&=&\frac{1}{2}x-\frac{1}{2}.
\end{eqnarray*}
The estimate of $D_{\rm \chi^2}(\boldsymbol{\theta},\boldsymbol{\theta}_0) $ is given by
\begin{eqnarray*}
\widehat{D}_{\rm \chi^2}(\boldsymbol{\theta},\boldsymbol{\theta}_0)&=&\sup_{\boldsymbol{\alpha}\in \mathbf{\Theta}}\left\{
\int\left(\frac{{\rm d}\mathbb{P}_{\boldsymbol{\theta}}}{{\rm d}\mathbb{P}_{\boldsymbol{\alpha}}}-1
\right){\rm d}\mathbb{P}_{\boldsymbol{\theta}}-\frac{1}{2}\int\left(\left(\frac{{\rm d}\mathbb{P}_{\boldsymbol{\theta}}}{{\rm d}\mathbb{P}_{\boldsymbol{\alpha}}}\right)^2-1
\right){\rm d}\mathbb{P}_n\right\}
\end{eqnarray*}
and the estimate of the parameter $\boldsymbol{\theta}_0$, with escort parameter $\boldsymbol{\theta}$, is defined by
\begin{equation*}
\widehat{\boldsymbol{\alpha}}_{\chi^{2}}(\boldsymbol{\theta}) := \arg\sup_{\boldsymbol{\alpha}\in \mathbf{\Theta}}\left\{
\int\left(\frac{{\rm d}\mathbb{P}_{\boldsymbol{\theta}}}{{\rm d}\mathbb{P}_{\boldsymbol{\alpha}}}-1
\right){\rm d}\mathbb{P}_{\boldsymbol{\theta}}-\frac{1}{2}\int\left(\left(\frac{{\rm d}\mathbb{P}_{\boldsymbol{\theta}}}{{\rm d}\mathbb{P}_{\boldsymbol{\alpha}}}\right)^2-1
\right){\rm d}\mathbb{P}_n
\right\}.
\end{equation*}
\item An other example is the Hellinger divergence
\begin{eqnarray*}
  \phi(x)&=&2(\sqrt{x}-1)^2\\
\phi^\prime(x)&=&2-\frac{1}{\sqrt{x}}\\
x  \phi^\prime(x)-\phi(x)&=&2\sqrt{x}-2.
\end{eqnarray*}
The estimate of $D_{\rm H}(\boldsymbol{\theta},\boldsymbol{\theta}_0)$ is given by
\begin{equation*}
\widehat{D}_{\rm H}(\boldsymbol{\theta},\boldsymbol{\theta}_0)=\sup_{\boldsymbol{\alpha}\in \mathbf{\Theta}}\left\{\int\left(2-2\sqrt{\frac{{\rm d}\mathbb{P}_{\boldsymbol{\alpha}}}{{\rm d}\mathbb{P}_{\boldsymbol{\theta}}}}\right){\rm d}\mathbb{P}_{\boldsymbol{\theta}}
-\int2\left(\sqrt{\frac{{\rm d}\mathbb{P}_{\boldsymbol{\theta}}}{{\rm d}\mathbb{P}_{\boldsymbol{\alpha}}}}-1
\right){\rm d}\mathbb{P}_n\right\}
\end{equation*}
and the estimate of the parameter $\boldsymbol{\theta}_0$, with escort parameter $\boldsymbol{\theta}$, is defined by
\begin{equation*}
\widehat{\boldsymbol{\alpha}}_{\rm H}(\boldsymbol{\theta}) := \arg\sup_{\boldsymbol{\alpha}\in \mathbf{\Theta}}\left\{\int\left(2-2\sqrt{\frac{{\rm d}\mathbb{P}_{\boldsymbol{\alpha}}}{{\rm d}\mathbb{P}_{\boldsymbol{\theta}}}}\right){\rm d}\mathbb{P}_{\boldsymbol{\theta}}
-\int2\left(\sqrt{\frac{{\rm d}\mathbb{P}_{\boldsymbol{\theta}}}{{\rm d}\mathbb{P}_{\boldsymbol{\alpha}}}}-1
\right){\rm d}\mathbb{P}_n
\right\}.
\end{equation*}

\item All the above examples are particular cases of the so-called
``\emph{power divergences}'',  which are defined through the class of convex real valued
functions, for $\gamma$ in $\mathbb{R}\backslash\left\{0,1\right\}$,
\begin{equation*}\label{gamma convexfunctions}
x\in\mathbb{R}_{+}^*\rightarrow
\varphi_{\gamma}(x):=\frac{x^{\gamma }-\gamma x+\gamma -1}{\gamma
(\gamma -1)}.
\end{equation*}
 The estimate of $D_{\gamma}(\boldsymbol{\theta},\boldsymbol{\theta}_0)$ is given by
\begin{eqnarray}\label{estimatorgammadiver}
\nonumber \widehat{D}_{\gamma}(\boldsymbol{\theta},\boldsymbol{\theta}_0)&=&\sup_{\boldsymbol{\alpha}\in \mathbf{\Theta}}\left\{\int\frac{1}{\gamma-1}\left(\left(\frac{{\rm d}\mathbb{P}_{\boldsymbol{\theta}}}{{\rm d}\mathbb{P}
_{\boldsymbol{\alpha}}}\right)^{\gamma-1}-1
\right){\rm d}\mathbb{P}_{\boldsymbol{\theta}}\right.\\&&\left.-\int\frac{1}{\gamma}\left(\left(\frac{{\rm d}\mathbb{P}_{\boldsymbol{\theta}}}{{\rm d}\mathbb{P}_{\boldsymbol{\alpha}}}\right)^\gamma-1
\right){\rm d}\mathbb{P}_n\right\}
\end{eqnarray}
and the parameter estimate  is defined by
\begin{eqnarray}
\lefteqn{\widehat{\boldsymbol{\alpha}}_{\gamma}(\boldsymbol{\theta})}\\&&:=- \arg\sup_{\boldsymbol{\alpha}\in \mathbf{\Theta}}
\left\{\int\frac{1}{\gamma-1}\left(\left(\frac{{\rm d}\mathbb{P}_{\boldsymbol{\theta}}}{{\rm d}\mathbb{P}_{\boldsymbol{\alpha}}}\right)^{\gamma-1}-1
\right){\rm d}\mathbb{P}_{\boldsymbol{\theta}}-\int\frac{1}{\gamma}\left(\left(\frac{{\rm d}\mathbb{P}_{\boldsymbol{\theta}}}{{\rm d}\mathbb{P}_{\boldsymbol{\alpha}}}\right)^\gamma-1
\right){\rm d}\mathbb{P}_n\right\}.\nonumber
\end{eqnarray}
\end{itemize}

\begin{remark} 
The computation of the estimate $\widehat{\boldsymbol{\alpha}}_\phi(\boldsymbol{\theta})$ requires calculus of the integral in the formula (\ref{Definition-h}). This integral can be explicitly  calculated for the most standard parametric models. 
Below,  we give a closed-form expression for Normal, log-Normal, Exponential, Gamma, Weilbull and Pareto density models. Hence, the computation of $\widehat{\boldsymbol{\alpha}}_\phi(\boldsymbol{\theta})$ can be performed by any standard non linear optimization code.
Unfortunately, the explicit formula of $\widehat{\boldsymbol{\alpha}}_\phi(\boldsymbol{\theta})$, generally, can not be derived, which also is the case for the ML method.  In practical problems, to obtain the estimate $\widehat{\boldsymbol{\alpha}}_\phi(\boldsymbol{\theta})$,
one can use the  Newton-Raphson algorithm taking as initial point the escort parameter $\boldsymbol{\theta}$. 
This algorithm, is a powerful technique for solving equations numerically, performs well since the the objective functions 
$\boldsymbol{\alpha}\in\mathbf{ \Theta}\mapsto \mathbb{P}_{\boldsymbol{\theta}_{0}}h(\boldsymbol{\theta},\boldsymbol{\alpha})$
are concave and the estimated parameter is unique for  functions  $\boldsymbol{\alpha}\in\mathbf{ \Theta}\mapsto \mathbb{P}_nh(\boldsymbol{\theta},\boldsymbol{\alpha})$, for instance, refer to \cite[Remark 3.5]{BroniatowskiKeziou2009}. 
\end{remark}

\subsection{Example of normal density} Consider the case of power divergences and the normal model
$$
\left\{N\left(\boldsymbol{\theta},\boldsymbol{\sigma}^{2}\right): (\boldsymbol{\theta},\boldsymbol{\sigma}^{2})\in\mathbf{ \Theta}=\mathbb{R}\times\mathbb{R}_{+}^{*}\right\}.
$$  
Set
$$
p_{\boldsymbol{\theta},\boldsymbol{\sigma}}(x)=\frac{1}{\boldsymbol{\sigma}\sqrt{2\pi}}\exp\left\{-\frac{1}{2}\left(\frac{x-\boldsymbol{\theta}}{\boldsymbol{\sigma}}\right)^{2}\right\}.
$$
Simple calculus gives, for $\gamma$  in $\mathbb{R}\backslash\{0,1\} $, 
\begin{eqnarray*}
\lefteqn{\frac{1}{\gamma-1}\int\left(\frac{{\rm d}\mathbb{P}_{\boldsymbol{\theta},\boldsymbol{\sigma}_{1}}(x)}{{\rm d}\mathbb{P}_{\boldsymbol{\alpha},\boldsymbol{\sigma}_{2}}(x)}\right)^{\gamma-1}~{\rm d}\mathbb{P}_{\boldsymbol{\theta},\boldsymbol{\sigma}_{1}}(x)dx}\\&&=
\frac{1}{\gamma-1}\frac{\boldsymbol{\sigma}_{1}^{-(\gamma-1)}\boldsymbol{\sigma}_{2}^{\gamma}}{\sqrt{\gamma\boldsymbol{\sigma}_{2}^{2}-(\gamma-1)\boldsymbol{\sigma}_{1}^{2}}}\exp\left\{\frac{\gamma(\gamma-1)(\boldsymbol{\theta}-\boldsymbol{\alpha})^{2}}{2(\gamma\boldsymbol{\sigma}_{2}^{2}-(\gamma-1)\boldsymbol{\sigma}_{1}^{2})}\right\}.
\end{eqnarray*}
This yields to
\begin{eqnarray*}
\lefteqn{\widehat{D}_{\gamma}((\boldsymbol{\theta},\boldsymbol{\sigma}_1),(\boldsymbol{\theta}_{0},\boldsymbol{\sigma}_0))}\\&=&\sup_{\boldsymbol{\alpha},\boldsymbol{\sigma}_{2}}\left\{\frac{1}{\gamma-1}\frac{\boldsymbol{\sigma}_{1}^{-(\gamma-1)}\boldsymbol{\sigma}_{2}^{\gamma}}{\sqrt{\gamma\boldsymbol{\sigma}_{2}^{2}-(\gamma-1)\boldsymbol{\sigma}_{1}^{2}}}\exp\left\{\frac{\gamma(\gamma-1)(\boldsymbol{\theta}-\boldsymbol{\alpha})^{2}}{2(\gamma\boldsymbol{\sigma}_{2}^{2}-(\gamma-1)\boldsymbol{\sigma}_{1}^{2})}\right\}\right.\\&&-\left.
\frac{1}{\gamma n}\sum_{i=1}^{n}\left(\frac{\boldsymbol{\sigma}_{2}}{\boldsymbol{\sigma}_{1}}\right)^{\gamma}\exp\left\{-\frac{\gamma}{2}\left(\left(\frac{\mathbf{ X}_{i}-\boldsymbol{\theta}}{\boldsymbol{\sigma}_{1}}\right)^{2}-\left(\frac{\mathbf{ X}_{i}-\boldsymbol{\alpha}}{\boldsymbol{\sigma}_{2}}\right)^{2}\right)
\right\} -\frac{1}{\gamma(\gamma-1)}\right\}.
\end{eqnarray*}
In the particular case, $\mathbb{P}_{\boldsymbol{\theta}}\equiv\mathcal{N}(\boldsymbol{\theta},1)$,
it follows that, for $\gamma\in\mathbb{R}\setminus\left\{0,1\right\}$,
\begin{eqnarray*}
\widehat{D}_{\gamma}\left(\boldsymbol{\theta},\boldsymbol{\theta}_{0}\right)&:=&\sup_{\boldsymbol{\alpha}}\int h\left(\boldsymbol{\theta},\boldsymbol{\alpha}\right){\rm{d}}\mathbb{P}_{n}\\
&=&\sup_{\boldsymbol{\alpha}}\left\{\frac{1}{\gamma-1}\exp\left\{\frac{\gamma(\gamma-1)(\boldsymbol{\theta}-\boldsymbol{\alpha})^2}{2}\right\}\right.\\
&&\left.-\frac{1}{\gamma n}\sum_{i=1}^n\exp\left\{-\frac{\gamma}{2}(\boldsymbol{\theta}-\boldsymbol{\alpha})(\boldsymbol{\theta}+\boldsymbol{\alpha}-2\mathbf{ X}_i)\right\}-\frac{1}{\gamma(\gamma-1)}\right\}. \end{eqnarray*}
For $\gamma=0$,
\begin{eqnarray*}
\widehat{D}_{\rm KL_m}\left(\boldsymbol{\theta},\boldsymbol{\theta}_{0}\right)&:=&\sup_{\boldsymbol{\alpha}}\int h\left(\boldsymbol{\theta},\boldsymbol{\alpha}\right){\rm{d}}\mathbb{P}_{n}\\
&=&\sup_{\boldsymbol{\alpha}}\left\{\frac{1}{2n}\sum_{i=1}^n\left(\boldsymbol{\theta}-\boldsymbol{\alpha}\right)\left(\boldsymbol{\theta}+\boldsymbol{\alpha}-2\mathbf{ X}_i\right)\right\},
\end{eqnarray*}
which lead to the maximum likelihood estimate independently upon $\boldsymbol{\theta}$.

\noindent For $\gamma=1$,
\begin{eqnarray*}
\widehat{D}_{\rm KL}\left(\boldsymbol{\theta},\boldsymbol{\theta}_{0}\right)&:=&\sup_{\boldsymbol{\alpha}}\int h\left(\boldsymbol{\theta},\boldsymbol{\alpha}\right){\rm{d}}\mathbb{P}_{n}\\
&=&\sup_{\boldsymbol{\alpha}}\left\{-\frac{1}{2}\left(\boldsymbol{\theta}-\boldsymbol{\alpha}\right)^2-\frac{1}{n}\sum_{i=1}^n\exp\left\{-\frac{1}{2}\left(\boldsymbol{\theta}-\boldsymbol{\alpha}\right)\left(\boldsymbol{\theta}+\boldsymbol{\alpha}-2\mathbf{ X}_i\right)\right\}+1\right\}.
\end{eqnarray*}

\subsection{Example of log-normal density} Consider the case of power divergences and the log-normal model
$$
\left\{p_{\boldsymbol{\theta},\boldsymbol{\sigma}}(x)=\frac{1}{x\boldsymbol{\sigma}\sqrt{2\pi}}\exp\left\{-\frac{1}{2}\left(\frac{\log(x)-\boldsymbol{\theta}}{\boldsymbol{\sigma}}\right)^{2}\right\}: (\boldsymbol{\theta},\boldsymbol{\sigma}^{2})\in\mathbf{ \Theta}=\mathbb{R}\times\mathbb{R}_{+}^{*}, x>0\right\}.
$$
Simple calculus gives, for $\gamma$  in $\mathbb{R}\backslash\{0,1\} $, 
\begin{eqnarray*}
\lefteqn{\frac{1}{\gamma-1}\int\left(\frac{{\rm d}\mathbb{P}_{\boldsymbol{\theta},\boldsymbol{\sigma}_{1}}(x)}{{\rm d}\mathbb{P}_{\boldsymbol{\alpha},\boldsymbol{\sigma}_{2}}(x)}\right)^{\gamma-1}~{\rm d}\mathbb{P}_{\boldsymbol{\theta},\boldsymbol{\sigma}_{1}}(x)dx}\\&&=
\frac{1}{\gamma-1}\frac{\boldsymbol{\sigma}_{1}^{-(\gamma-1)}\boldsymbol{\sigma}_{2}^{\gamma}}{\sqrt{\gamma\boldsymbol{\sigma}_{2}^{2}-(\gamma-1)\boldsymbol{\sigma}_{1}^{2}}}\exp\left\{\frac{\gamma(\gamma-1)(\boldsymbol{\theta}-\boldsymbol{\alpha})^{2}}{2(\gamma\boldsymbol{\sigma}_{2}^{2}-(\gamma-1)\boldsymbol{\sigma}_{1}^{2})}\right\}.
\end{eqnarray*}
This yields to
\begin{eqnarray*}
\lefteqn{\widehat{D}_{\gamma}((\boldsymbol{\theta},\boldsymbol{\sigma}_1),(\boldsymbol{\theta}_{0},\boldsymbol{\sigma}_0))}\\&=&\sup_{\boldsymbol{\alpha},\boldsymbol{\sigma}_{2}}\left\{\frac{1}{\gamma-1}\frac{\boldsymbol{\sigma}_{1}^{-(\gamma-1)}\boldsymbol{\sigma}_{2}^{\gamma}}{\sqrt{\gamma\boldsymbol{\sigma}_{2}^{2}-(\gamma-1)\boldsymbol{\sigma}_{1}^{2}}}\exp\left\{\frac{\gamma(\gamma-1)(\boldsymbol{\theta}-\boldsymbol{\alpha})^{2}}{2(\gamma\boldsymbol{\sigma}_{2}^{2}-(\gamma-1)\boldsymbol{\sigma}_{1}^{2})}\right\}\right.\\&&-\left.
\frac{1}{\gamma n}\sum_{i=1}^{n}\left(\frac{\boldsymbol{\sigma}_{2}}{\boldsymbol{\sigma}_{1}}\right)^{\gamma}\exp\left\{-\frac{\gamma}{2}\left(\left(\frac{\log(\mathbf{ X}_{i})-\boldsymbol{\theta}}{\boldsymbol{\sigma}_{1}}\right)^{2}-\left(\frac{\log(\mathbf{ X}_{i})-\boldsymbol{\alpha}}{\boldsymbol{\sigma}_{2}}\right)^{2}\right)
\right\}\right.\\ &&-\left.\frac{1}{\gamma(\gamma-1)}\right\}.
\end{eqnarray*}
 \subsection{Example of exponential density}
  Consider the case of power divergences and the exponential model
$$
\left\{p_{\boldsymbol{\theta}}(x)=\boldsymbol{\theta}\exp(-\boldsymbol{\theta} x): \boldsymbol{\theta}\in \Theta=\mathbb{R}_{+}^{*}\right\}.
$$  
We have, for $\gamma$  in $\mathbb{R}\backslash\{0,1\} $, 
 \begin{equation}\label{exampleexp1}
 \frac{1}{\gamma-1}\int\left(\frac{{\rm d}\mathbb{P}_{\boldsymbol{\theta}}(x)}{{\rm d}\mathbb{P}_{\boldsymbol{\alpha}}(x)}\right)^{\gamma-1}~{\rm d}\mathbb{P}_{\boldsymbol{\theta}}(x)dx=\left(\frac{\boldsymbol{\theta}}{\boldsymbol{\alpha}}\right)^{(\gamma-1)}\left(\frac{\boldsymbol{\theta}}{\boldsymbol{\theta}\gamma(\gamma-1)-\boldsymbol{\alpha}(\gamma-1)^{2}}\right).
 \end{equation}
 Then using this last equality, one finds 
\begin{eqnarray*}
\widehat{D}_{\gamma}(\boldsymbol{\theta},\boldsymbol{\theta}_{0})&=&\sup_{\boldsymbol{\alpha}}\left\{\left(\frac{\boldsymbol{\theta}}{\boldsymbol{\alpha}}\right)^{(\gamma-1)}\left(\frac{\boldsymbol{\theta}}{\boldsymbol{\theta}\gamma(\gamma-1)-\boldsymbol{\alpha}(\gamma-1)^{2}}\right)\right.\\&&-\left.
\frac{1}{\gamma n}\sum_{i=1}^{n}\left(\frac{\boldsymbol{\theta}}{\boldsymbol{\alpha}}\right)^{\gamma}\exp\left\{-\gamma\left(\left(\boldsymbol{\theta}\mathbf{ X}_{i}\right)-\left(\boldsymbol{\alpha}\mathbf{ X}_{i}\right)\right)
\right\} -\frac{1}{\gamma(\gamma-1)}\right\}.
\end{eqnarray*}
In more  general case,  we may consider the gamma density combined with the power divergence. The Gamma model is defined by 
$$
\left\{p_{\boldsymbol{\theta}}(x;k):=\boldsymbol{\theta}^{k}x^{k-1}\frac{\exp{(-x\boldsymbol{\theta}})}{\Gamma(k)}:k,\boldsymbol{\theta}\geq 0 \right\},
$$
where $\Gamma(\cdot)$ is the Gamma function 
$$
\Gamma(k):=\int_{0}^{\infty}x^{k-1}\exp(-x)dx.
$$
Simple calculus gives, for $\gamma$  in $\mathbb{R}\backslash\{0,1\} $,  
\begin{equation*}
 \frac{1}{\gamma-1}\int\left(\frac{{\rm d}\mathbb{P}_{\boldsymbol{\theta};k}(x)}{{\rm d}\mathbb{P}_{\boldsymbol{\alpha};k}(x)}\right)^{\gamma-1}~{\rm d}\mathbb{P}_{\boldsymbol{\theta};k}(x)dx=\left(\frac{\boldsymbol{\theta}}{\boldsymbol{\alpha}}\right)^{k(\gamma-1)}\left(\frac{\boldsymbol{\theta}}{\boldsymbol{\theta}\gamma-\boldsymbol{\alpha}(\gamma-1)}\right)^{k}\frac{1}{\gamma-1}, 
 \end{equation*}
 which implies that
\begin{eqnarray*}
\widehat{D}_{\gamma}(\boldsymbol{\theta},\boldsymbol{\theta}_{0})&=&\sup_{\boldsymbol{\alpha}}\left\{\left(\frac{\boldsymbol{\theta}}{\boldsymbol{\alpha}}\right)^{k(\gamma-1)}\left(\frac{\boldsymbol{\theta}}{\boldsymbol{\theta}\gamma-\boldsymbol{\alpha}(\gamma-1)}\right)^{k}\frac{1}{\gamma-1}\right.\\&&-\left.
\frac{1}{\gamma n}\sum_{i=1}^{n}\left(\frac{\boldsymbol{\theta}}{\boldsymbol{\alpha}}\right)^{k\gamma}\exp\left\{-\gamma\left(\left(\boldsymbol{\theta}\mathbf{ X}_{i}\right)-\left(\boldsymbol{\alpha}\mathbf{ X}_{i}\right)\right)
\right\} -\frac{1}{\gamma(\gamma-1)}\right\}.
\end{eqnarray*}
\subsection{Example of Weibull density}
Consider the case of power divergences and the Weibull density model, with the assumption   that $k\in \mathbb{R}_{+}^{*}$ is known and $\boldsymbol{\theta}$ is the parameter of interest to be estimated, recall that 
$$
\left\{p_{\boldsymbol{\theta}}(x)=\frac{k}{\boldsymbol{\theta}}\left(\frac{x}{\boldsymbol{\theta}}\right)^{k-1}\exp\left(-\left(\frac{x}{\boldsymbol{\theta}}\right)^{k}\right): \boldsymbol{\theta}\in \mathbf{ \Theta}=\mathbb{R}_{+}^{*}, x\geq 0\right\}.
$$
Routine algebra gives, for $\gamma$  in $\mathbb{R}\backslash\{0,1\} $,  
\begin{equation*}
 \frac{1}{\gamma-1}\int\left(\frac{{\rm d}\mathbb{P}_{\boldsymbol{\theta};k}(x)}{{\rm d}\mathbb{P}_{\boldsymbol{\alpha};k}(x)}\right)^{\gamma-1}~{\rm d}\mathbb{P}_{\boldsymbol{\theta};k}(x)dx=\left(\frac{\boldsymbol{\alpha}}{\boldsymbol{\theta}}\right)^{k(\gamma-1)}\left(\frac{1}{\gamma-\left(\frac{\boldsymbol{\theta}}{\boldsymbol{\alpha}}\right)^{k}(\gamma-1)}\right)\frac{1}{\gamma-1}, 
 \end{equation*}
 which implies that
\begin{eqnarray*}
\widehat{D}_{\gamma}(\boldsymbol{\theta},\boldsymbol{\theta}_{0})&=&\sup_{\boldsymbol{\alpha}}\left\{\left(\frac{\boldsymbol{\alpha}}{\boldsymbol{\theta}}\right)^{k(\gamma-1)}\left(\frac{1}{\gamma-\left(\frac{\boldsymbol{\theta}}{\boldsymbol{\alpha}}\right)^{k}(\gamma-1)}\right)\frac{1}{\gamma-1}\right.\\&&-\left.
\frac{1}{\gamma n}\sum_{i=1}^{n}\left(\frac{\boldsymbol{\alpha}}{\boldsymbol{\theta}}\right)^{k\gamma}\exp\left\{-\gamma\left(\left(\frac{\mathbf{ X}_{i}}{\boldsymbol{\theta}}\right)^{k}-\left(\frac{\mathbf{ X}_{i}}{\boldsymbol{\alpha}}\right)^{k}\right)
\right\} -\frac{1}{\gamma(\gamma-1)}\right\}.
\end{eqnarray*}

\subsection{Example of the Pareto density} Consider the case of power divergences and the Pareto density
$$
\left\{p_{\boldsymbol{\theta}}(x):=\frac{\boldsymbol{\theta}}{x^{\boldsymbol{\theta}+1}}:~x>1;~\boldsymbol{\theta}\in\mathbb{R}_{+}^{*}\right\}.
$$  
Simple calculus gives, for $\gamma$  in $\mathbb{R}\backslash\{0,1\} $, 
 \begin{equation}\label{examplepareto1}
 \frac{1}{\gamma-1}\int\left(\frac{{\rm d}\mathbb{P}_{\boldsymbol{\theta}}(x)}{{\rm d}\mathbb{P}_{\boldsymbol{\alpha}}(x)}\right)^{\gamma-1}~{\rm d}\mathbb{P}_{\boldsymbol{\theta}}(x)~{\rm d}x=\left(\frac{\boldsymbol{\theta}}{\boldsymbol{\alpha}}\right)^{(\gamma-1)}\left(\frac{\boldsymbol{\theta}}{\boldsymbol{\theta}\gamma(\gamma-1)-\boldsymbol{\alpha}(\gamma-1)^{2}}\right).
 \end{equation}
As before, using this last equality, one finds 
 \begin{eqnarray*}
 \widehat{D}_{\gamma}(\boldsymbol{\theta},\boldsymbol{\theta}_{0})&=&\sup_{\boldsymbol{\alpha}}\left\{\left(\frac{\boldsymbol{\theta}}{\boldsymbol{\alpha}}\right)^{(\gamma-1)}\left(\frac{\boldsymbol{\theta}}{\boldsymbol{\theta}\gamma(\gamma-1)-\boldsymbol{\alpha}(\gamma-1)^{2}}\right)\right.\\&&-\left.
 \frac{1}{\gamma n}\sum_{i=1}^{n}\left(\frac{\boldsymbol{\theta}}{\boldsymbol{\alpha}}\right)^{\gamma}\mathbf{ X}_{i}^{\left\{-\gamma\left(\boldsymbol{\theta}-\boldsymbol{\alpha}\right)
 \right\}} -\frac{1}{\gamma(\gamma-1)}\right\}.
 \end{eqnarray*}
For $\gamma=0$,
\begin{eqnarray*}
\widehat{D}_{\rm KL_m}\left(\boldsymbol{\theta},\boldsymbol{\theta}_{0}\right)&:=&\sup_{\boldsymbol{\alpha}}\int h\left(\boldsymbol{\theta},\boldsymbol{\alpha}\right){\rm{d}}\mathbb{P}_{n}\\
&=&\sup_{\boldsymbol{\alpha}}\left\{-\frac{1}{n}\sum_{i=1}^n\left\{\log\left(\frac{\boldsymbol{\theta}}{\boldsymbol{\alpha}}\right)-\left(\boldsymbol{\theta}-\boldsymbol{\alpha}\right)\log\left(\mathbf{ X}_{i}\right)\right\}\right\},
\end{eqnarray*}
which lead to the maximum likelihood estimate, given by $$\left(\frac{1}{n}\sum_{i=1}^n\log\left(\mathbf{ X}_{i}\right)\right)^{-1},$$ independently upon $\boldsymbol{\theta}$.

\begin{remark}
The choice of divergence, i.e., the statistical criterion,  depends crutially on the problem at hand. For example, the $\chi^2$-divergence among various divergences in the nonstandard problem (e.g., boundary problem estiamtion) is more appropriate. The idea is to include the parameter domain
$\mathbf{ \Theta}$ into an enlarged space, say $\mathbf{ \Theta}_e$, in order to
render the boundary value an interior point of the new parameter space,
$\mathbf{ \Theta}_e$. Indeed, 
Kullback-Leibler, modified Kullback-Leibler, modified
$\chi^2$, and Hellinger  divergences  are infinite when ${\rm d}\mathbb{Q}/{\rm d}\mathbb{P}$ takes negative values on non negligible  (with respect to $\mathbb{P}$) subset of the support of $\mathbb{P}$, since the corresponding
$\phi(\cdot)$ is infinite on $(-\infty,0)$, when $\boldsymbol{\theta}$ belongs  to $\mathbf{ \Theta}_e\backslash\mathbf{ \Theta}$. This problem does not hold in the case of $\chi^2$-divergence, in fact, the corresponding $\phi(\cdot)$ is finite on $\mathbb{R}$, for more details refer to \cite{bouzebda-keziou2008,BouzebdaKeziou2010b,BouzebdaKeziou2010a}, consult also   \cite{BroniatowskiKeziou2009} and \cite{Broniatowski2006} for related matter. 
It is well known that when the underlying model is misspecified or when the data are contaminated the maximum likelihood or other classical parametric methods may be severely affected and lead to very poor results.  Therefore, robust methods, which automatically circumvent the   contamination  effects and model misspecification, can be used to provide a compromise between efficient classical parametric methods and the semi-parametric approach provided they are reasonably efficient at the model, this problem has been investigated in \cite{Basu1998,Basu2006}.
In \cite{BouzebdaKeziou2010b,BouzebdaKeziou2010a},
simulation
results show that the choice of $\chi^2$-divergence has good
properties in terms of efficiency-robustness.
 We mention that  some progress has been made on automatic data-based selection of the tuning parameter $\alpha>0$, appearing in  formula (1) of  \cite{Basu2006},  the interested reader is referred to  \cite{Hong2001} and \cite{Jones2005}. It is mentioned in  \cite{tsukahara2005}, where semiparametric minimum distance estimators are considered, that  the MLE or inversion-type estimators involve solving a nonlinear equation which depends on some initial value. The second difficulty is that the objective function is not convex in $\boldsymbol{\theta}$, in general, which give the situation of multiple roots. Thus in general, ``\emph{good}'' consistent initial estimate are necessary and the D$\phi$DE should serve that purpose.
\end{remark}
\section{Random right censoring}\label{censored}
\noindent Let $T=T_{1}, \ldots,T_{n}$ be i.i.d. survival times with
continuous survival function $1-F_{\boldsymbol{\theta}_{0}}(\cdot)=1-\mathbb{P}_{\boldsymbol{\theta}_{0}}(T\leq \cdot)$  and
$C_{1},\ldots,C_{n}$ be independent censoring times with d.f.
$G(\cdot)$. In the censoring set-up, we observe only the pair
$Y_{i}= \min(T_{i},C_{i})$ and $\delta_{i}=\mathds{1}\{T_{i}\leq
C_{i}\}$, where $\mathds{1}\{\cdot\}$ is the indicator function of the event $\{\cdot\}$, which designs whether an observation has been censored
or not. Let $(Y_{1},\delta_{1}),\ldots,(Y_{n},\delta_{n})$ denote
the observed data points and $$t(1)< t(2)<\cdots< t(k)$$ be the $k$
distinct death times. Now define the death set and risk set as
follows, for $j=1,\ldots, k$,
$$D(j):=\{i:y_{i}=t(j),\delta_{i}=1\}$$
and
$$R(j):=\{i:y_{i}\geq t(j)\}.$$
The \cite{kaplanmeier1958}'s estimator of $1-F_{\theta_{0}}(\cdot)$, denoted here by $1-\widehat{F}_n(\cdot)$,  may be written as follows
$$
1-\widehat{F}_n(t):=\prod_{j=1}^{k}\left(1-\frac{\sum_{q\in D(j)}1}{\sum_{q\in R(j)}1}\right)^{\mathds{1}\{T_{(j)}\leq t\}}.
$$
One may define a generally
exchangeable weighted bootstrap scheme for the Kaplan-Meier estimator
and related functionals as follows, cf. \cite[p. 1598]{Lancelot1997},
$$
1-\widehat{F}_{n}^{*}(t):=\prod_{j=1}^{k}\left(1-\frac{\sum_{q\in D(j)}
W_{nq}}{\sum_{q\in R(j)}W_{nq}}\right)^{\mathds{1}\{T(j)\leq t\}}.
$$
Let $\psi$ be $F_{\boldsymbol{\theta}_{0}}$-integrable and put
$$
\Psi_{n}:=\int\psi(u)\mathrm{d}\widehat{\mathbb{P}}_{n}^{*}(u)=\sum_{j=1}^{k}\Upsilon_{jn}\psi(T_{(j)}),
$$
where
$$
\Upsilon_{jn}:=\left(\frac{\sum_{q\in D(j)}W_{nq}}{\sum_{q\in R(j)}W_{nq}}\right)
\prod_{k=1}^{j-1}\left(\frac{\sum_{q\in D(k)}W_{nq}}{\sum_{q\in R(k)}W_{nq}}\right).
$$
Note that we have used the following identity.
Let $a_i$, $i=1,\ldots,k$, $b_i$, $i=1,\ldots,k$, be real numbers
\begin{equation*}\label{iden}
\prod_{i=1}^{k}a_i-\prod_{i=1}^{k}b_i=\sum_{i=1}^k(a_i-b_i)\prod_{j=1}^{i-1}b_j\prod_{h=1+i}^{k}a_h.
\end{equation*}
In the similar way, we define a more appropriate representation, that will be used in the sequel, as follows
$$
\Psi_{n}=\int\psi(u)\mathrm{d}\widehat{\mathbb{P}}_{n}^{*}(u)=\sum_{j=1}^{n}\pi_{jn}\psi(Y_{j:n}),
$$
where, for $1 \leq j \leq n$,
$$
\pi_{jn}:=\delta_{j:n}\left(\frac{\sum_{q\in D(j)}W_{nq}}{\sum_{q\in R(j)}W_{nq}}\right)\prod_{k=1}^{j-1}\left(\frac{\sum_{q\in D(k)}W_{nq}}{\sum_{q\in R(k)}W_{nq}}\right)^{\delta_{k:n}}.
$$
Here, $Y_{1:n}\leq \cdots\leq Y_{n:n}$ are ordered $Y$-values and
$\delta_{i:n}$ denotes the concomitant associated with $Y_{i:n}$.
Hence we may write
\begin{equation}\label{definitionPbootstap}
\widehat{\mathbb{P}}^{*}_{n}:=\sum_{j=1}^{n}\pi_{jn}\delta_{Y_{i:n}}
\end{equation}
For the right censoring situation, the bootstrap D$\phi$DE's, is defined by
replacing $\mathbb{P}_n$ in (\ref{dualestimator}) by
$\widehat{\mathbb{P}}^{*}_{n}$, that is
\begin{equation}\label{dualestimatorcensored}
\widehat{\boldsymbol{\alpha}}_{n}(\boldsymbol{\theta}):=\arg\sup_{\boldsymbol{\alpha}\in \mathbf{\Theta}}\int
h(\boldsymbol{\theta},\boldsymbol{\alpha})\mathrm{d}\widehat{\mathbb{P}}^{*}_{n},\;\;\boldsymbol{\theta}\in \mathbf{\Theta}.
\end{equation}
\noindent The corresponding estimating equation for the unknown parameter is then given by
\begin{equation}\label{estimatingequation}
\int\frac{\partial}{\partial\boldsymbol{\alpha}}h(\boldsymbol{\theta},\boldsymbol{\alpha})\mathrm{d}\widehat{\mathbb{P}}_{n}^{*}=0,
\end{equation}
where we recall that 
\begin{equation*}
h(\boldsymbol{\theta},\boldsymbol{\alpha},x):=\int \phi ^{\prime }\left( \frac{\mathrm{d}\mathbb{P}_{\boldsymbol{\theta} }}{\mathrm{d}\mathbb{P}_{\boldsymbol{\alpha}
}}\right) ~\mathrm{d}\mathbb{P}_{\boldsymbol{\theta} }-\left[ \frac{\mathrm{d}\mathbb{P}_{\boldsymbol{\theta} }(x)}{\mathrm{d}\mathbb{P}_{\boldsymbol{\alpha}
}(x)}\phi ^{\prime }\left( \frac{\mathrm{d}\mathbb{P}_{\boldsymbol{\theta} }(x)}{\mathrm{d}\mathbb{P}_{\boldsymbol{\alpha} }(x)}\right)
-\phi \left( \frac{\mathrm{d}\mathbb{P}_{\boldsymbol{\theta} } (x)}{\mathrm{d}\mathbb{P}_{\boldsymbol{\alpha} }(x)}\right) \right].
\end{equation*}
\noindent Formula (\ref{dualestimatorcensored}) defines a family of $M$-estimator
 for censored data. 
\noindent In the case of the power divergences family (\ref{powerdivergence}), it follows that from (\ref{estimatorgammadiver})
\begin{equation*}
\int
h(\boldsymbol{\theta},\boldsymbol{\alpha})\mathrm{d}\widehat{\mathbb{P}}_{n}=\frac{1}{\gamma-1}\int\left(\frac{\mathrm{d}\mathbb{P}_{\boldsymbol{\theta}}}{\mathrm{d}\mathbb{P}_{\boldsymbol{\alpha}}}\right)^{\gamma-1}~{\rm{d}}\mathbb{P}_{\boldsymbol{\theta}}-\frac{1}{\gamma}
\int\left[\left(\frac{\mathrm{d}\mathbb{P}_{\boldsymbol{\theta}}}{\mathrm{d}\mathbb{P}_{\boldsymbol{\alpha}}}\right)^{\gamma}-1\right]~{\rm{d}}\widehat{\mathbb{P}}_{n}-\frac{1}{\gamma-1},
\end{equation*}
where $$
\widehat{\mathbb{P}}_{n}:=\sum_{j=1}^{n}\omega_{jn}\delta_{Y_{i:n}},
$$
and, for $1 \leq j \leq n$,
\begin{equation*}
 \omega_{jn}=\frac{\delta_{j:n}}{n-j+1}\prod_{i=1}^{j-1}\left[\frac{n-i}{n-i+1}\right]^{\delta_{i:n}}.
\end{equation*}

\noindent Consider the lifetime distribution to be the one parameter exponential $\exp{(\boldsymbol{\theta})}$ with density
$\boldsymbol{\theta} e^{-\boldsymbol{\theta} x},~x\geq 0$.  Following \cite{Stute1995}, the {K}aplan-{M}eier integral $\int
h(\boldsymbol{\theta},\boldsymbol{\alpha})\mathrm{d}\widehat{\mathbb{P}}_{n}$  may be written as
\begin{equation*}
    \sum_{j=1}^n\omega_{jn}h(\boldsymbol{\theta},\boldsymbol{\alpha},Y_{j:n}).
\end{equation*}
The MLE of $\boldsymbol{\theta}_0$ is given by
\begin{equation}\label{MLE}
   \widehat{ \boldsymbol{\theta}}_{n,MLE}=\frac{\sum_{j=1}^n\delta_{j}}{\sum_{j=1}^nY_{j}},
\end{equation}
and the approximate MLE (AMLE) of \cite{Oakes1986} is defined by
\begin{equation}\label{AMLE}
   \widehat{ \boldsymbol{\theta}}_{n,AMLE}=\frac{\sum_{j=1}^n\delta_{j}}{\sum_{j=1}^n\omega_{jn}Y_{j:n}}.
\end{equation}
We infer from (\ref{exampleexp1}), that, for $\gamma\in\mathbb{R}\setminus\left\{0,1\right\}$,
\begin{eqnarray*}
 \int h(\boldsymbol{\theta},\boldsymbol{\alpha})\mathrm{d}\widehat{\mathbb{P}}_{n}&=&\frac{\boldsymbol{\theta}^\gamma\boldsymbol{\alpha}^{1-\gamma}}{(\gamma-1)\left[\gamma\boldsymbol{\theta}+(1-\gamma)\boldsymbol{\alpha}\right]}\\
&&-\frac{1}{\gamma}\sum_{j=1}^n\omega_{jn}\left[\left(\frac{\boldsymbol{\theta}}{\boldsymbol{\alpha}}\right)^\gamma\exp\left\{-\gamma(\boldsymbol{\theta}-\boldsymbol{\alpha})Y_{j:n}\right\}-1\right].
\end{eqnarray*}
For $\gamma=0$,
\begin{equation*}
 \int
h(\boldsymbol{\theta},\boldsymbol{\alpha})\mathrm{d}\widehat{\mathbb{P}}_{n}=\sum_{j=1}^n\omega_{jn}\left[(\boldsymbol{\theta}-\boldsymbol{\alpha})Y_{j:n}-\log\left(\frac{\boldsymbol{\theta}}{\boldsymbol{\alpha}}\right)\right].
\end{equation*}
Observe that this divergence leads to the AMLE, independently upon the value of $\boldsymbol{\theta}$.

\noindent For $\gamma=1$,
\begin{equation*}
 \int
h(\boldsymbol{\theta},\boldsymbol{\alpha})\mathrm{d}\widehat{\mathbb{P}}_{n}=\log\left(\frac{\boldsymbol{\theta}}{\boldsymbol{\alpha}}\right)-\frac{(\boldsymbol{\theta}-\boldsymbol{\alpha})}{\boldsymbol{\theta}}-\sum_{i=1}^n\omega_{jn}\left[ \frac{\boldsymbol{\theta}}{\boldsymbol{\alpha}}\exp\left(-(\boldsymbol{\theta}-\boldsymbol{\alpha})Y_{j:n}\right)-1
\right].
\end{equation*}
For more details  about dual
$\phi$-divergence estimators  in right censoring we refer to \cite{Cherfi2011b},
 we leave this study open for future research. We mention that the bootstrapped estimators, in this framework,  are obtained by replacing the weights  $\omega_{jn}$ by $\pi_{jn}$ in the preceding formulas.

\FloatBarrier

\section{Simulations}\label{simulation}

\noindent In this section, series of experiments were conducted in order to examine the performance of the proposed random weighted bootstrap procedure of the D$\phi$DE's, defined in (\ref{Bdualestimator}). We provide numerical illustrations regarding the mean squared error (MSE) and the coverage
probabilities. The computing program codes were implemented in \texttt{R}.

\noindent The values of $\gamma$ are chosen to be $-1,~0,~0.5,~1,~2$, which corresponds, as indicated above, to the well known standard divergences: $\chi_{m}^{2}$-divergence, $KL_{m}$, the {H}ellinger distance, $KL$  and the $\chi^{2}$-divergence respectively.  The samples of sizes considered in our simulations are $25,~50,~75,~100,~150,~200$ 
and the estimates, D$\phi$DE's $\widehat{\boldsymbol{\alpha}}_{\phi}(\boldsymbol{\theta})$, are obtained from $500$ independent runs.
The  value of escort parameter $\theta$ is taken to be the MLE, which, under the model, is a consistent estimate of $\theta_{0}$, and the limit distribution of the D$\phi$DE $\widehat{\boldsymbol{\alpha}}_{\phi}(\boldsymbol{\theta}_{0})$, in this case, has variance which indeed coincides with the MLE, for more details on this subject, we refer to \cite[Theorem 2.2, (1) (b)]{Keziou2003}, as it is mentioned in Section \ref{choiceoftheescort}. The bootstrap weights are chosen to be  $$(W_{n1},\ldots,W_{nn})\sim ~~\mbox{ Dirichlet}(n;1,\ldots,1).$$

\begin{figure}[!h]
\begin{center}
\centerline{
\includegraphics[width=8cm]{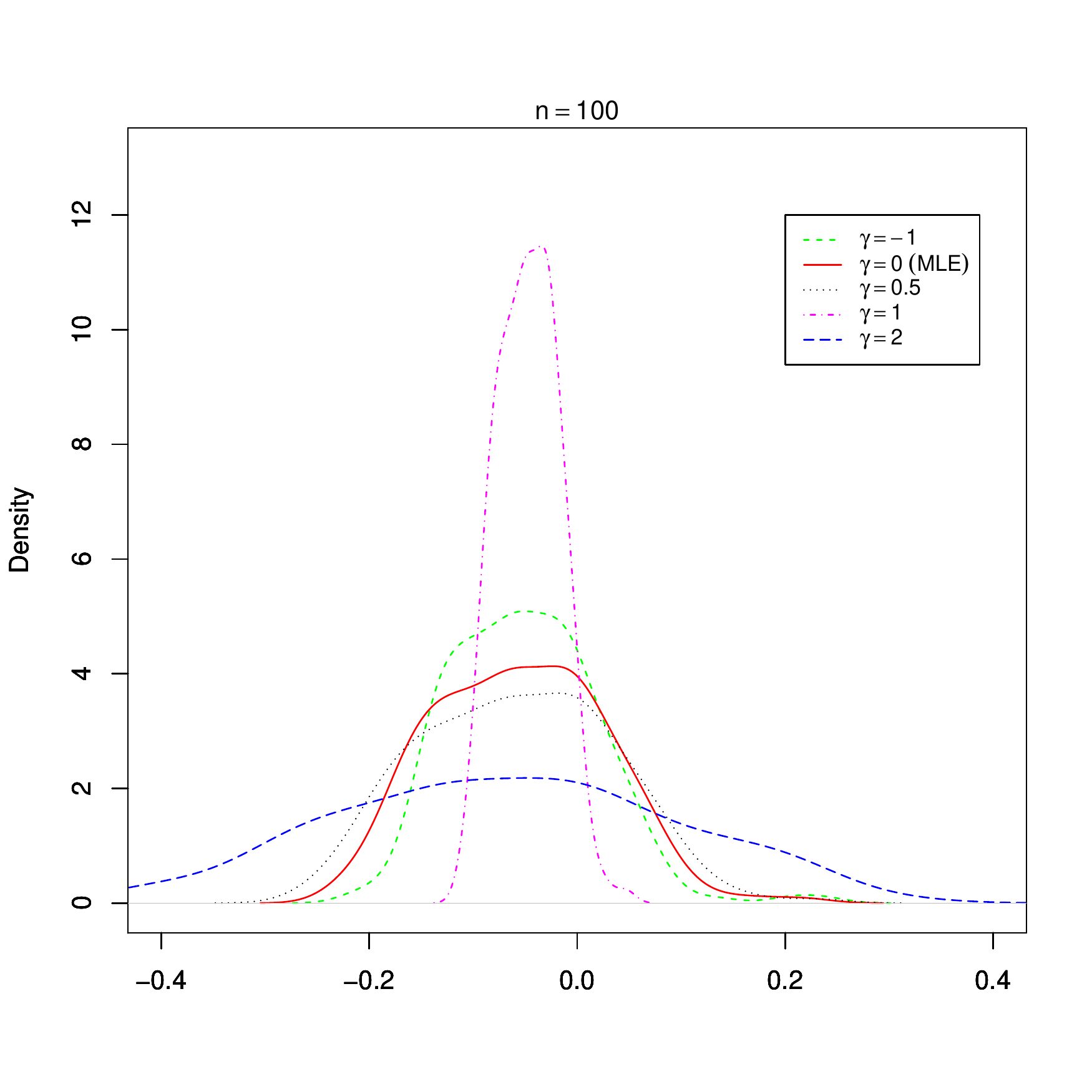}
\includegraphics[width=8cm]{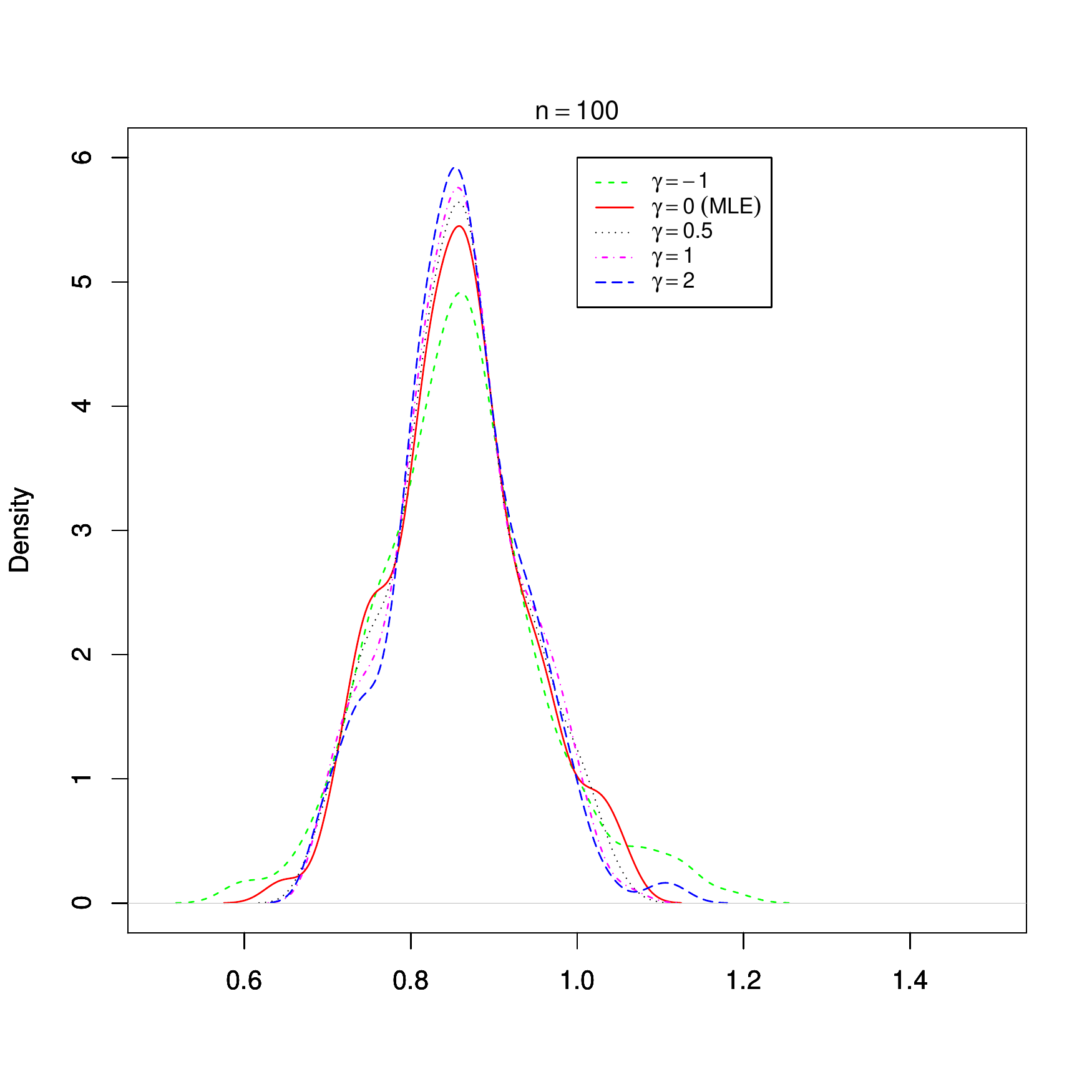}}
\end{center}
\caption{Densities of the estimates.} \label{FigSim1}
\end{figure}

 \noindent In Figure \ref{FigSim1}, we plot the densities of the different estimates, it shows that the proposed estimators perform 
reasonably well.

\noindent Tables \ref{TabSim1500}  and \ref{TabSim1} provide the MSE of various estimates under the Normal model $N(\theta_0=0,1)$. Here, we mention that the $KL$ based estimator ($\gamma=1$) is more efficient than the others competitors. 

\noindent Tables \ref{TabSim3500} and  \ref{TabSim3} provide the MSE of various estimates under the Exponential model $\exp(\theta_0=1)$. As
expected,  the MLE produces most efficient estimators. A close look at the results of the simulations show that the D$\phi$DE's
perform well under the model. For large sample size $n=200$, the estimator based on the {H}ellinger distance is equivalent to that of the MLE. Indeed in terms of empirical MSE the D$\phi$DE with $\gamma=0.5$ produces the same MSE as the MLE, while the performance of the other estimators is comparable.
\begin{table}[ht]
\caption{MSE of the estimates for the Normal distribution, B=500}
\begin{center}
\begin{tabular}{rrrrrrr}
  \hline
& $n=25$ & $n=50$ & $n=75$ & $n=100$ & $n=150$ & $n=200$ \\
 \cline{2-7}
  \cline{2-7}
   $\gamma$&&&&&&\\
-1 & 0.0687 & 0.0419 & 0.0288 & 0.0210 & 0.0135 & 0.0107 \\
 0 & 0.0647 & 0.0373 & 0.0255 & 0.0192 & 0.0127 & 0.0101 \\
 0.5 & 0.0668 & 0.0379 & 0.0257 & 0.0194 & 0.0128 & 0.0101 \\
 1 & 0.0419 & 0.0217 & 0.0143 & 0.0108 & 0.0070 & 0.0057 \\
 2 & 0.0931 & 0.0514 & 0.0331 & 0.0238 & 0.0148 & 0.0112 \\
  \hline
\end{tabular}
\end{center} \label{TabSim1500}
\end{table}

\begin{table}[ht]
\caption{MSE of the estimates for the Normal distribution, B=1000}
\begin{center}
\begin{tabular}{rrrrrrr}
 \hline
& $n=25$ & $n=50$ & $n=75$ & $n=100$ & $n=150$ & $n=200$ \\
  \cline{2-7}
   \cline{2-7}
  $\gamma$&&&&&&\\
-1 & 0.0716 & 0.0432 & 0.0285 & 0.0224 & 0.0147 & 0.0099 \\
 0 & 0.0670 & 0.0385 & 0.0255 & 0.0202 & 0.0136 & 0.0093 \\
 0.5 & 0.0684 & 0.0391 & 0.0258 & 0.0203 & 0.0137 & 0.0093 \\
 1 & 0.0441 & 0.0230 & 0.0143 & 0.0116 & 0.0078 & 0.0049 \\
 2 & 0.0900 & 0.0522 & 0.0335 & 0.0246 & 0.0156 & 0.0103 \\
  \hline
\end{tabular}
\end{center}\label{TabSim1}
\end{table}

\begin{table}[ht]
\caption{MSE of the estimates for the Exponential distribution, B=500}
\begin{center}
\begin{tabular}{rrrrrrr}
  \hline
& $n=25$ & $n=50$ & $n=75$ & $n=100$ & $n=150$ & $n=200$ \\
 \cline{2-7}
  \cline{2-7}
   $\gamma$&&&&&&\\
-1 & 0.0729 & 0.0435 & 0.0313 & 0.0215 & 0.0146 & 0.0117 \\
 0 & 0.0708 & 0.0405 & 0.0280 & 0.0195 & 0.0131 & 0.0104 \\
 0.5 & 0.0727 & 0.0415 & 0.0282 & 0.0197 & 0.0131 & 0.0105 \\
 1 & 0.0786 & 0.0446 & 0.0296 & 0.0207 & 0.0136 & 0.0108 \\
 2 & 0.1109 & 0.0664 & 0.0424 & 0.0289 & 0.0178 & 0.0132 \\
  \hline
\end{tabular}
\end{center} \label{TabSim3500}
\end{table}

\begin{table}[ht]
\caption{MSE of the estimates for the Exponential distribution, B=1000}
\begin{center}
\begin{tabular}{rrrrrrr}
 \hline
& $n=25$ & $n=50$ & $n=75$ & $n=100$ & $n=150$ & $n=200$ \\
 \cline{2-7}
  \cline{2-7}
   $\gamma$&&&&&&\\
-1 & 0.0670 & 0.0444 & 0.0295 & 0.0243 & 0.0146 & 0.0111 \\
 0 & 0.0659 & 0.0417 & 0.0269 & 0.0216 & 0.0133 & 0.0102 \\
 0.5 & 0.0677 & 0.0427 & 0.0272 & 0.0216 & 0.0135 & 0.0102 \\
 1 & 0.0735 & 0.0458 & 0.0287 & 0.0225 & 0.0140 & 0.0106 \\
 2 & 0.1074 & 0.0697 & 0.0429 & 0.0306 & 0.0183 & 0.0133 \\
  \hline
\end{tabular}
\end{center}\label{TabSim3}
\end{table}

\noindent Tables \ref{TabSim2500}, \ref{TabSim2}, \ref{TabSim4500} and \ref{TabSim4}, provide the empirical coverage
probabilities of the corresponding $0.95$ weighted bootstrap confidence
intervals based on $B=500, 1000$ weighted bootstrap estimators. 
Notice that the empirical coverage
probabilities as in any other inferential
context, the greater the sample size, the better. From the results reported in these tables, we find that for large values of the sample size $n$, the empirical coverage probabilities are all close to  the nominal level. One can see that the D$\phi$DE with $\gamma=2$ has the best empirical coverage
probability which is near the assigned nominal level.
\begin{table}[ht]
\caption{Empirical coverage
probabilities for the Normal distribution, B=500 }
\begin{center}
\begin{tabular}{rrrrrrr}
  \hline
& $n=25$ & $n=50$ & $n=75$ & $n=100$ & $n=150$ & $n=200$ \\
 \cline{2-7}
  \cline{2-7}
   $\gamma$&&&&&&\\
-1 & 0.88 & 0.91 & 0.93 & 0.92 & 0.95 & 0.92 \\
 0 & 0.91 & 0.92 & 0.94 & 0.94 & 0.94 & 0.93 \\
 0.5 & 0.94 & 0.94 & 0.94 & 0.96 & 0.94 & 0.93 \\
 1 & 0.44 & 0.47 & 0.54 & 0.46 & 0.48 & 0.51 \\
 2 & 0.97 & 0.97 & 0.96 & 0.97 & 0.95 & 0.95 \\
  \hline
\end{tabular}
\end{center} \label{TabSim2500}
\end{table}

\begin{table}[ht]
\caption{Empirical coverage
probabilities for the Normal distribution, B=1000 }
\begin{center}
\begin{tabular}{rrrrrrr}
 \hline
& $n=25$ & $n=50$ & $n=75$ & $n=100$ & $n=150$ & $n=200$ \\
 \cline{2-7}
  \cline{2-7}
   $\gamma$&&&&&&\\
-1 & 0.87 & 0.90 & 0.93 & 0.92 & 0.93 & 0.96 \\
 0 & 0.91 & 0.94 & 0.94 & 0.93 & 0.94 & 0.96 \\
 0.5 & 0.93 & 0.93 & 0.95 & 0.93 & 0.94 & 0.96 \\
 1 & 0.46 & 0.45 & 0.48 & 0.46 & 0.45 & 0.50 \\
 2 & 0.96 & 0.97 & 0.96 & 0.95 & 0.96 & 0.96 \\
  \hline
\end{tabular}
\end{center}\label{TabSim2}
\end{table}

\begin{table}[ht]
 \caption{Empirical coverage
probabilities for the Exponential distribution, B=500}
\begin{center}
\begin{tabular}{rrrrrrr}
 \hline
& $n=25$ & $n=50$ & $n=75$ & $n=100$ & $n=150$ & $n=200$ \\
 \cline{2-7}
  \cline{2-7}
   $\gamma$&&&&&&\\
-1 & 0.67 & 0.83 & 0.87 & 0.91 & 0.93 & 0.92 \\
 0 & 0.73 & 0.87 & 0.91 & 0.93 & 0.96 & 0.93 \\
 0.5 & 0.76 & 0.88 & 0.91 & 0.94 & 0.96 & 0.93 \\
 1 & 0.76 & 0.88 & 0.90 & 0.95 & 0.97 & 0.93 \\
 2 & 0.76 & 0.89 & 0.91 & 0.96 & 0.96 & 0.94 \\
  \hline
\end{tabular}
\end{center}\label{TabSim4500}
\end{table}
\begin{table}[ht]
 \caption{Empirical coverage
probabilities for the Exponential distribution, B=1000}
\begin{center}
\begin{tabular}{rrrrrrr}
 \hline
& $n=25$ & $n=50$ & $n=75$ & $n=100$ & $n=150$ & $n=200$ \\
 \cline{2-7}
  \cline{2-7}
   $\gamma$&&&&&&\\
-1 & 0.70 & 0.79 & 0.90 & 0.91 & 0.92 & 0.91 \\
 0 & 0.76 & 0.84 & 0.91 & 0.92 & 0.93 & 0.92 \\
 0.5 & 0.78 & 0.85 & 0.93 & 0.94 & 0.94 & 0.93 \\
 1 & 0.78 & 0.87 & 0.94 & 0.94 & 0.95 & 0.94 \\
 2 & 0.78 & 0.88 & 0.95 & 0.95 & 0.96 & 0.95 \\
  \hline
\end{tabular}
\end{center}\label{TabSim4}
\end{table}

\subsection{Right censoring case}

\noindent  This subsection presents some simulations for right censoring case discussed in \S \ref{censored}.  A sample is generated from $\exp(1)$ and an exponential censoring scheme is used, the censoring distribution  is taken to be $\exp(1/9)$, that the proportion of censoring is $10\%$. To study the robustness properties of our estimators $20\%$ of the observations are contaminated by $\exp(5)$. The D$\phi$DE's $\widehat{\boldsymbol{\alpha}}_{\phi}(\boldsymbol{\theta})$ are calculated for samples of sizes $25,~50,~100,~150$ and the hole procedure is repeated $500$ times.
\begin{table}[ht]
\caption{MSE of the estimates for the Exponential distribution under right censoring}
\begin{center}
\begin{tabular}{rrrrr}
  \hline
& $n=25$ & $n=50$ &$n=100$ & $n=150$\\
 \cline{2-5}
  \cline{2-5}
   $\gamma$&&&&\\
-1 & 0.1088 & 0.0877 & 0.0706 & 0.0563\\
0 & 0.1060 & 0.0843 & 0.0679 & 0.0538\\
0.5 & 0.1080 & 0.0860 & 0.0689 & 0.0544\\
1 & 0.1150 & 0.0914 & 0.0724 & 0.0567\\
2& 0.1535 & 0.1276 & 0.1019 & 0.0787\\
   \hline
\end{tabular}
\end{center}\label{TabSimCens1}
\end{table}
We can see from Table \ref{TabSimCens1} that the D$\phi$DE's
perform well under the model in term of MSE, and are an attractive alternative to the AMLE.

\begin{table}[ht]
 \caption{Empirical coverage
probabilities for the Exponential distribution under right censoring}
\begin{center}
\begin{tabular}{rrrrr}
 \hline
& $n=25$ & $n=50$ & $n=100$ & $n=150$ \\
 \cline{2-5}
  \cline{2-5}
   $\gamma$&&&&\\
-1 & 0.55 & 0.63 & 0.63 & 0.64 \\
0 & 0.59 & 0.66 & 0.64 & 0.64\\
0.5 & 0.61 & 0.66 & 0.64 & 0.65\\
1 & 0.63 & 0.67 & 0.66 & 0.66\\
2 & 0.64 & 0.70 &  0.68 & 0.67\\
\hline
\end{tabular}
\end{center}\label{TabSimCens2}
\end{table}

Table \ref{TabSimCens2}  shows the variation
in coverage of nominal $95\%$ asymptotic confidence intervals according to
the sample size. There clearly
is under coverage of the confidence intervals, the D$\phi$DE's  have
poor coverage probabilities due to the censoring effect. However for small and moderate sized
samples the  D$\phi$DE's associated to $\gamma=2$  outperforms the AMLE. 

\begin{table}[ht]
\caption{MSE of the estimates for the Exponential distribution under right censoring, $20\%$ of contamination}
\begin{center}
\begin{tabular}{rrrrr}
  \hline
& $n=25$ & $n=50$ &$n=100$ & $n=150$\\
 \cline{2-5}
  \cline{2-5}
   $\gamma$&&&&\\
-1 & 0.1448 & 0.1510 & 0.1561 & 0.1591\\
0 & 0.1482 & 0.1436 & 0.1409 & 0.1405\\
0.5 & 0.1457 & 0.1402 & 0.1360 & 0.1342\\
1 & 0.1462 & 0.1389 & 0.1332 & 0.1300\\
2& 0.1572 & 0.1442 & 0.1338 & 0.1266\\
   \hline
\end{tabular}
\end{center}\label{TabSimCens3}
\end{table}

\begin{table}[ht]
 \caption{Empirical coverage
probabilities for the Exponential distribution under right censoring, $20\%$ of contamination}
\begin{center}
\begin{tabular}{rrrrr}
 \hline
& $n=25$ & $n=50$ & $n=100$ & $n=150$ \\
 \cline{2-5}
  \cline{2-5}
   $\gamma$&&&&\\
-1 & 0.44 & 0.49 & 0.54 & 0.57 \\
0 & 0.46 & 0.49 & 0.53 & 0.57\\
0.5 & 0.46 & 0.49 & 0.53 & 0.57\\
1 & 0.45 & 0.49 & 0.53 & 0.57\\
2 & 0.45 & 0.49 &  0.52 & 0.53\\
\hline
\end{tabular}
\end{center}\label{TabSimCens4}
\end{table}

Under contamination the performances of our estimators decrease considerably. Such findings are evidences for the need of more adequate procedures for right censored data.

\begin{remark}
In order to extract methodological recommendations for the use of an appropriate divergence, it will be interesting to conduct an extensive Monte Carlo experiments for several divergences or investigate theoretically the problem of the choice
of the divergence which leads to an ``\emph{optimal}'' (in some sense)
estimate in terms of efficiency and robustness, which would go well beyond the scope of the present paper. An other challenging task is how to choose the bootstrap weights for a given divergence in order to obtain, for example, an efficient estimator.
\end{remark}

\section{Appendix}\label{proof}
\noindent This section is devoted to the proofs of our results. The previously defined notation continues to be used below.
\subsection{Proof of Theorem \ref{asythm-a}}
Proceeding as \cite{vanderVaartWellner1996} in their proof  of  the Argmax theorem, i.e., Corollary  3.2.3, it is straightforward to show the
consistency  of the bootstrapped estimates $\widehat{\boldsymbol{\alpha}}^{\ast}_{\phi}(\boldsymbol{\theta})$.
\begin{flushright}
$\Box$
\end{flushright}
\begin{remark}
Note that the proof techniques of Theorem \ref{asythm-b} are largely inspired from that of \cite{ChengHuang2010} and changes have been made in order to
adapt them to our purpose.
\end{remark}

\subsection{Proof of Theorem \ref{asythm-b}}
Keep in mind the following definitions
$$\mathbb{G}_n:=\sqrt{n}(\mathbb{P}_n-\mathbb{P}_{\boldsymbol{\theta}_0})$$ and
$$\mathbb{G}_n^{\ast}:=\sqrt{n}(\mathbb{P}_n^{\ast}-\mathbb{P}_n).$$
In view of the fact that $\displaystyle{\mathbb{P}_{\boldsymbol{\theta}_0}\frac{\partial}{\partial\boldsymbol{\alpha}}h(\boldsymbol{\theta},\boldsymbol{\theta}_0)=0}$, then a little calculation shows that
\begin{eqnarray*}
\lefteqn{\mathbb{G}_{n}^{\ast}\frac{\partial}{\partial\boldsymbol{\alpha}}h(\boldsymbol{\theta},\boldsymbol{\theta}_0)+\mathbb{G}_{n}\frac{\partial}{\partial\boldsymbol{\alpha}}h(\boldsymbol{\theta},\boldsymbol{\theta}_0)}\\&&+\sqrt{n}
\mathbb{P}_{\boldsymbol{\theta}_0}\left[\frac{\partial}{\partial\boldsymbol{\alpha}}h(\boldsymbol{\theta},\widehat{\boldsymbol{\alpha}}^{\ast}_{\phi}(\boldsymbol{\theta}))
-\frac{\partial}{\partial\boldsymbol{\alpha}}h(\boldsymbol{\theta},\boldsymbol{\theta}_0)\right]\\
&=&\mathbb{G}_{n}^{\ast}\left[\frac{\partial}{\partial\boldsymbol{\alpha}}h(\boldsymbol{\theta},\boldsymbol{\theta}_0)-\frac{\partial}{\partial\boldsymbol{\alpha}}h(\boldsymbol{\theta},\widehat{\boldsymbol{\alpha}}^{\ast}_{\phi}(\boldsymbol{\theta}))\right]\\&&+\mathbb{G}_{n}\left[\frac{\partial}{\partial\boldsymbol{\alpha}}h(\boldsymbol{\theta},\boldsymbol{\theta}_0)-
\frac{\partial}{\partial\boldsymbol{\alpha}}h(\boldsymbol{\theta},\widehat{\boldsymbol{\alpha}}^{\ast}_{\phi}(\boldsymbol{\theta}))\right]\\&&+\sqrt{n}\mathbb{P}
_{n}^{\ast}\frac{\partial}{\partial\boldsymbol{\alpha}}h(\boldsymbol{\theta},\widehat{\boldsymbol{\alpha}}^{\ast}_{\phi}(\boldsymbol{\theta})).
\end{eqnarray*}
Consequently,  we have following inequality
\begin{eqnarray}\label{equarat}
\nonumber\lefteqn{\left\|\sqrt{n}\mathbb{P}_{\boldsymbol{\theta}_0}\left[\frac{\partial}{\partial\boldsymbol{\alpha}}h(\boldsymbol{\theta},\widehat{\boldsymbol{\alpha}}^{\ast}_{\phi}(\boldsymbol{\theta}))-\frac{\partial}{\partial\boldsymbol{\alpha}}h(\boldsymbol{\theta},\boldsymbol{\theta}_0)\right]\right\|}\\
&\leq&\left\|
\mathbb{G}_{n}^{\ast}\frac{\partial}{\partial\boldsymbol{\alpha}}h(\boldsymbol{\theta},\boldsymbol{\theta}_0)\right\|+\left\|\mathbb{G}_{n}\frac{\partial}{\partial\boldsymbol{\alpha}}h(\boldsymbol{\theta},\boldsymbol{\theta}_0)\right\|
\nonumber\\&&+\left\|\mathbb{G}_{n}^{\ast}\left[\frac{\partial}{\partial\boldsymbol{\alpha}}h(\boldsymbol{\theta},\widehat{\boldsymbol{\alpha}}^{\ast}_{\phi}(\boldsymbol{\theta}))-\frac{\partial}{\partial\boldsymbol{\alpha}}h(\boldsymbol{\theta},\boldsymbol{\theta}_0)\right]\right\|\nonumber\\\nonumber
&&+\left\|\mathbb{G}_{n}\left[\frac{\partial}{\partial\boldsymbol{\alpha}}h(\boldsymbol{\theta},\widehat{\boldsymbol{\alpha}}^{\ast}_{\phi}(\boldsymbol{\theta}))-\frac{\partial}{\partial\boldsymbol{\alpha}}h(\boldsymbol{\theta},\boldsymbol{\theta}_0)\right]\right\|\\\nonumber&&+
\left\|\sqrt{n}\mathbb{P}_{n}^{\ast}\frac{\partial}{\partial\boldsymbol{\alpha}}h(\boldsymbol{\theta},\widehat{\boldsymbol{\alpha}}^{\ast}_{\phi}(\boldsymbol{\theta}))\right\|\\
&:=& G_{1}+G_{2}+G_{3}+G_{4}+G_{5}.
\end{eqnarray}
According to Theorem 2.2 in \cite{PraestgaardWellner1993}, under
condition {\rm{(A.4)}}, we have $G_{1}=O_{\mathbb{P}_{W}}^{o}(1)$
in $\mathbb{P}_{\boldsymbol{\theta}_0}$-probability. In view of the CLT,
we have $G_{2}=O_{\mathbb{P}_{\boldsymbol{\theta}_0}}(1)$.

\noindent By applying a Taylor series expansion, we have
\begin{equation}\label{eqaa}
\mathbb{G}_{n}^{\ast}\left[\frac{\partial}{\partial\boldsymbol{\alpha}}h(\boldsymbol{\theta},\widehat{\boldsymbol{\alpha}}^{\ast}_{\phi}(\boldsymbol{\theta}))-\frac{\partial}{\partial\boldsymbol{\alpha}}h(\boldsymbol{\theta},\boldsymbol{\theta}_0)\right]=\left(\widehat{\boldsymbol{\alpha}}^{\ast}_{\phi}(\boldsymbol{\theta})-\boldsymbol{\theta}_0\right)^\top\mathbb{G}_{n}^{\ast}\frac{\partial^2}{\partial\boldsymbol{\alpha}^2}h(\boldsymbol{\theta},\overline{\boldsymbol{\alpha}}),
\end{equation}
where $\overline{\boldsymbol{\alpha}}$ is between
$\widehat{\boldsymbol{\alpha}}^{\ast}_{\phi}(\boldsymbol{\theta})$ and
$\boldsymbol{\theta}_0$. By condition {\rm{(A.5)}} and  Theorem
2.2 in \cite{PraestgaardWellner1993}, we conclude that the  right term in
(\ref{eqaa})  is of order
$O_{\mathbb{P}_{W}}^{o}\left(\|\widehat{\boldsymbol{\alpha}}^{\ast}_{\phi}(\boldsymbol{\theta})-\boldsymbol{\theta}_{0}\|\right)$
in $\mathbb{P}_{\boldsymbol{\theta}_0}$-probability. The fact that
$\widehat{\boldsymbol{\alpha}}^{\ast}_{\phi}(\boldsymbol{\theta})$ is assumed to be
consistent, then, we have $G_{3}=o_{\mathbb{P}_{W}}^{o}(1)$ in
$\mathbb{P}_{\boldsymbol{\theta}_0}$-probability.
\noindent An
analogous argument yields
$$\displaystyle{\mathbb{G}_{n}\left[\frac{\partial}{\partial\boldsymbol{\alpha}}h(\boldsymbol
{\theta},\widehat{\boldsymbol{\alpha}}^{\ast}_{\phi}(\boldsymbol{\theta}))-\frac{\partial}{\partial\boldsymbol{\alpha}}h(\boldsymbol{\theta},\boldsymbol{\theta}_0)\right]}$$
is of order
$O_{\mathbb{P}_{\boldsymbol{\theta}_0}}\left(\|\widehat{\boldsymbol{\alpha}}^{\ast}_{\phi}(\boldsymbol{\theta})-\boldsymbol{\theta}_{0}\|\right)$,
by the consistency of $\widehat{\boldsymbol{\alpha}}^{\ast}_{\phi}(\boldsymbol{\theta})$, we
have $G_{4}=o_{\mathbb{P}_{W}}^{o}(1)$ in
$\mathbb{P}_{\boldsymbol{\theta}_0}$-probability. Finally,
$G_{5}=0$ based on (\ref{solBdualestimator}). In summary,
(\ref{equarat}) can be rewritten as follows
\begin{eqnarray}
\left\|\sqrt{n}\mathbb{P}_{\boldsymbol{\theta}_0}(\frac{\partial}{\partial\boldsymbol{\alpha}}
h(\boldsymbol{\theta},\widehat{\boldsymbol{\alpha}}^{\ast}_{\phi}(\boldsymbol{\theta}))-\frac{\partial}{\partial\boldsymbol{\alpha}}
h(\boldsymbol{\theta},\boldsymbol{\theta}_0))\right\|\leq
O_{\mathbb{P}_{W}}^{o}(1)+O_{\mathbb{P}_{\boldsymbol{\theta}_0}}^{o}(1)\label{equarat2}
\end{eqnarray}
in $\mathbb{P}_{\boldsymbol{\theta}_0}$-probability. On the other
hand, by  a Taylor series expansion, we can write

\begin{equation}\label{Taylor2}
\mathbb{P}_{\boldsymbol{\theta}_0}\left[\frac{\partial}{\partial\boldsymbol{\alpha}}
h(\boldsymbol{\theta},\boldsymbol{\alpha})-\frac{\partial}{\partial\boldsymbol{\alpha}}
h(\boldsymbol{\theta},\boldsymbol{\theta}_0)\right]=-(\boldsymbol{\alpha}-\boldsymbol{\theta}_0)^\top S
+O\left(\|\boldsymbol{\alpha}-\boldsymbol{\theta}_{0}\|^2\right).
\end{equation}
Clearly it is straightforward to combine (\ref{Taylor2}) with (\ref{equarat2}), to infer the following
\begin{eqnarray}
\sqrt{n} \left\|S\|\widehat{\boldsymbol{\alpha}}^{\ast}_{\phi}(\boldsymbol{\theta})-\boldsymbol{\theta}_{0}\|\right\|\leq
O_{\mathbb{P}_{W}}^{o}(1)+O_{\mathbb{P}_{\boldsymbol{\theta}_0}}^{o}(1)+O_{\mathbb{P}_{W}}^{o}
\left(\sqrt{n}\|\widehat{\boldsymbol{\alpha}}^{\ast}_{\phi}(\boldsymbol{\theta})-\boldsymbol{\theta}_{0}\|^{2}\right)\label{rootnpf}
\end{eqnarray}
in $\mathbb{P}_{\boldsymbol{\theta}_0}$-probability. By
considering again the consistency of $\widehat{\boldsymbol{\alpha}}^{\ast}_{\phi}(\boldsymbol{\theta})$ and
condition {\rm{(A.3)}} and making use  (\ref{rootnpf}) to complete the proof of (\ref{bconratep}).

\noindent We next prove (\ref{bcons}). Introduce
\begin{eqnarray*}
H_1&:=&-\mathbb{G}_{n}^{\ast}\left[\frac{\partial}{\partial\boldsymbol{\alpha}}h(\boldsymbol{\theta},\widehat{\boldsymbol{\alpha}}^{\ast}_{\phi}(\boldsymbol{\theta}))-\frac{\partial}{\partial\boldsymbol{\alpha}}h(\boldsymbol{\theta},\boldsymbol{\theta}_0)\right],\\
H_2&:=&\mathbb{G}_{n}\left[\frac{\partial}{\partial\boldsymbol{\alpha}}h(\boldsymbol{\theta},\widehat{\boldsymbol{\alpha}}_{\phi}(\boldsymbol{\theta}))-\frac{\partial}{\partial\boldsymbol{\alpha}}h(\boldsymbol{\theta},\boldsymbol{\theta}_0)\right],\\
H_3&:=&-\mathbb{G}_{n}\left[\frac{\partial}{\partial\boldsymbol{\alpha}}h(\boldsymbol{\theta},\widehat{\boldsymbol{\alpha}}^{\ast}_{\phi}(\boldsymbol{\theta}))-\frac{\partial}{\partial\boldsymbol{\alpha}}h(\boldsymbol{\theta},\boldsymbol{\theta}_0)\right],\\
H_4&:=&\sqrt{n}\mathbb{P}_{n}^{\ast}\frac{\partial}{\partial\boldsymbol{\alpha}}h(\boldsymbol{\theta},\widehat{\boldsymbol{\alpha}}^{\ast}_{\phi}(\boldsymbol{\theta}))-\sqrt{n}\mathbb{P}_{n}
\frac{\partial}{\partial\boldsymbol{\alpha}}h(\boldsymbol{\theta},\widehat{\boldsymbol{\alpha}}_{\phi}(\boldsymbol{\theta})).
\end{eqnarray*}
By some algebra, we obtain
$$\sqrt{n}\mathbb{P}_{\boldsymbol{\theta}_0}\left(\frac{\partial}{\partial\boldsymbol{\alpha}}h(\boldsymbol{\theta},\widehat{\boldsymbol{\alpha}}^{\ast}_{\phi}(\boldsymbol{\theta}))-\frac{\partial}{\partial\boldsymbol{\alpha}}h(\boldsymbol{\theta},\widehat{\boldsymbol{\alpha}}_{\phi}(\boldsymbol{\theta}))\right)
+\mathbb{G}_{n}^{\ast}\frac{\partial}{\partial\boldsymbol{\alpha}}h(\boldsymbol{\theta},\boldsymbol{\theta}_0)=\sum_{j=1}^4 H_j.$$
Obviously, $H_1=O^{o}_{\mathbb{P}_{W}}(n^{-1/2})$ in $\mathbb{P}_{\boldsymbol{\theta}_0}$-probability and $H_2=O_{\mathbb{P}_{\boldsymbol{\theta}_0}}(n^{-1/2})$. We also know that the order of $H_3$ is $O^{o}_{\mathbb{P}_{W}}(n^{-1/2})$ in $\mathbb{P}_{\boldsymbol{\theta}_0}$-probability. Using (\ref{soldualestimator}) and (\ref{solBdualestimator}) we obtain that $H_4=0$.

\noindent Therefore, we have established
\begin{eqnarray}
\nonumber\sqrt{n}\mathbb{P}_{\boldsymbol{\theta}_0}\left[\frac{\partial}{\partial\boldsymbol{\alpha}}h(\boldsymbol{\theta},\widehat{\boldsymbol{\alpha}}^{\ast}_{\phi}(\boldsymbol{\theta}))-\frac{\partial}{\partial\boldsymbol{\alpha}}h(\boldsymbol{\theta},\widehat{\boldsymbol{\alpha}}_{\phi}(\boldsymbol{\theta}))\right]&=&-
\mathbb{G}_{n}^{\ast}\frac{\partial}{\partial\boldsymbol{\alpha}}h(\boldsymbol{\theta},\boldsymbol{\theta}_0)+
o_{\mathbb{P}_{\boldsymbol{\theta}_0}}(1)\\&&+o_{\mathbb{P}_{W}}^{o}(1) \label{inter5}
\end{eqnarray}
in $\mathbb{P}_{\boldsymbol{\theta}_0}$-probability. To analyze the left hand side of
(\ref{inter5}), we rewrite it as
$$\sqrt{n}\mathbb{P}_{\boldsymbol{\theta}_0}\left[\frac{\partial}{\partial\boldsymbol{\alpha}}h(\boldsymbol{\theta},\widehat{\boldsymbol{\alpha}}^{\ast}_{\phi}(\boldsymbol{\theta}))-\frac{\partial}{\partial\boldsymbol{\alpha}}h(\boldsymbol{\theta},\boldsymbol{\theta}_0)\right]-\sqrt{n}
\mathbb{P}_{\boldsymbol{\theta}_0}\left[\frac{\partial}{\partial\boldsymbol{\alpha}}h(\boldsymbol{\theta},\widehat{\boldsymbol{\alpha}}_{\phi}(\boldsymbol{\theta}))-\frac{\partial}{\partial\boldsymbol{\alpha}}h(\boldsymbol{\theta},\boldsymbol{\theta}_0)\right].$$
By a Taylor expansion, we obtain
\begin{eqnarray}\label{inter6}
\nonumber\lefteqn{\sqrt{n}S(\widehat{\boldsymbol{\alpha}}^{\ast}_{\phi}(\boldsymbol{\theta})-\widehat{\boldsymbol{\alpha}}_{\phi}(\boldsymbol{\theta}))}\\\nonumber&=&
\mathbb{G}_{n}^{\ast}\frac{\partial}{\partial\boldsymbol{\alpha}}h(\boldsymbol{\theta},\boldsymbol{\theta}_0)+
o_{\mathbb{P}_{\boldsymbol{\theta}_0}}(1)+o_{\mathbb{P}_{W}}^{o}(1)\\\nonumber&&+O_{\mathbb{P}_{\boldsymbol{\theta}_0}}(n^{-1/2})+O_{\mathbb{P}_{W}}^{o}(n^{-1/2})\\
&=&
\mathbb{G}_{n}^{\ast}\frac{\partial}{\partial\boldsymbol{\alpha}}h(\boldsymbol{\theta},\boldsymbol{\theta}_0)+
o_{\mathbb{P}_{\boldsymbol{\theta}_0}}(1)+o_{\mathbb{P}_{W}}^{o}(1)
\end{eqnarray}
in $\mathbb{P}_{\boldsymbol{\theta}_0}$-probability.
Keep in mind that, under condition {\rm{(A.3)}}, the matrix $S$ is nonsingular. Multiply both sides of (\ref{inter6}) by $S^{-1}$
to obtain (\ref{bcons}). An application of  \cite[Lemma 4.6]{PraestgaardWellner1993}, under the bootstrap weight conditions, thus implies (\ref{bconcor}).  Using  \cite[Theorem 3.2]{BroniatowskiKeziou2009} and  \cite[Lemma 2.11]{vanderVaart1998}, it easily follows that
\begin{equation}\label{convq}
\sup_{\mathbf{ x}\in\mathbb{R}^d}\left|\mathbb{P}_{\boldsymbol{\theta}_0}(\sqrt{n}(\widehat{\boldsymbol{\alpha}}_{\phi}(\boldsymbol{\theta})-\boldsymbol{\theta}_0)\leq
\mathbf{ x})-\mathbb{P}(N(0,\Sigma)\leq
\mathbf{ x})\right|=o_{\mathbb{P}_{\boldsymbol{\theta}_0}}(1).
\end{equation}
By combining (\ref{bconcor}) and (\ref{convq}), we readily obtain the desired conclusion (\ref{proconv}).
\begin{flushright}
$\Box$
\end{flushright}

\section*{Acknowledgements} 
\noindent  We are grateful to the referees, whose insightful comments help to improve an early draft of this article greatly. The authors are indebted to Amor 
 Keziou for careful reading and fruitful discussions on the subject.
 We would like to thank the Associate Editor for comments which helped in the completion of this work.

\def\ocirc#1{\ifmmode\setbox0=\hbox{$#1$}\dimen0=\ht0 \advance\dimen0
  by1pt\rlap{\hbox to\wd0{\hss\raise\dimen0
  \hbox{\hskip.2em$\scriptscriptstyle\circ$}\hss}}#1\else {\accent"17 #1}\fi}


\begin{thebibliography}{}
\bibitem[Barbe and Bertail(1995)]{Bertail95}
Barbe, P. and Bertail, P. (1995).
\newblock {\em The weighted bootstrap}, volume~98 of {\em Lecture Notes in
  Statistics}.
\newblock Springer-Verlag, New York.

\bibitem[Basu {\em et~al.}(1998)]{Basu1998}
Basu, A., Harris, I.~R., Hjort, N.~L., and Jones, M.~C. (1998).
\newblock Robust and efficient estimation by minimising a density power
  divergence.
\newblock {\em Biometrika}, {\bf 85}(3), 549--559.

\bibitem[Basu {\em et~al.}(2006)]{Basu2006}
Basu, S., Basu, A., and Jones, M.~C. (2006).
\newblock Robust and efficient parametric estimation for censored survival
  data.
\newblock {\em Ann. Inst. Statist. Math.}, {\bf 58}(2), 341--355.

\bibitem[Basu {\em et~al.}(2011)Basu, Shioya, and Park]{BasuShioyaPark2011}
Basu, A., Shioya, H., and Park, C. (2011).
\newblock {\em Statistical Inference: The Minimum Distance Approach\/}, volume
  120 of {\em Monographs on Statistics \& Applied Probability\/}.
\newblock Chapman \& Hall/CRC, Boca Raton, FL.

\bibitem[Beran(1984)]{Beran1984}
Beran, R. (1984).
\newblock Bootstrap methods in statistics.
\newblock {\em Jahresber. Deutsch. Math.-Verein.}, {\bf 86}(1), 14--30.

\bibitem[Beran and Millar(1986)]{BeranMillar1986}
Beran, R. and Millar, P.~W. (1986).
\newblock Confidence sets for a multivariate distribution.
\newblock {\em Ann. Statist.}, {\bf 14}(2), 431--443.

\bibitem[Beran {\em et~al.}(1987)]{BeranLeCamMillar1987}
Beran, R.~J., Le~Cam, L., and Millar, P.~W. (1987).
\newblock Convergence of stochastic empirical measures.
\newblock {\em J. Multivariate Anal.}, {\bf 23}(1), 159--168.

\bibitem[Bickel and Freedman(1981)]{BickelFreedman1981}
Bickel, P.~J. and Freedman, D.~A. (1981).
\newblock Some asymptotic theory for the bootstrap.
\newblock {\em Ann. Statist.}, {\bf 9}(6), 1196--1217.

\bibitem[Bouzebda and Keziou(2008)]{bouzebda-keziou2008}
Bouzebda, S. and Keziou, A. (2008).
\newblock A test of independence in some copula models.
\newblock {\em Math. Methods Statist.}, {\bf 17}(2), 123--137.



\bibitem[Bouzebda and Keziou(2010a)]{BouzebdaKeziou2010b}
Bouzebda, S. and Keziou, A. (2010a).
\newblock Estimation and tests of independence in copula models via
  divergences.
\newblock {\em Kybernetika}, {\bf 46}(1), 178-201.

\bibitem[Bouzebda and Keziou(2010b)]{BouzebdaKeziou2010a}
Bouzebda, S. and Keziou, A. (2010b).
\newblock A new test procedure of independence in copula models via
  {$\chi^2$}-divergence.
\newblock {\em Comm. Statist. Theory Methods}, {\bf 39}, 1--20.

\bibitem[Broniatowski(2011)]{Broniatowski2011}
Broniatowski, M. (2011).
\newblock Minimum divergence estimators, maximum likelihood and exponential families.
\newblock {\em Arxiv preprint:  arXiv:1108.0772}.

\bibitem[Broniatowski and Leorato(2006)]{Broniatowski2006}
Broniatowski, M. and Leorato, S. (2006).
\newblock An estimation method for the {N}eyman chi-square divergence with
  application to test of hypotheses.
\newblock {\em J. Multivariate Anal.}, {\bf 97}(6), 1409--1436.


\bibitem[Broniatowski and Keziou(2006)]{BroniatowskiKeziou2006}
Broniatowski, M. and Keziou, A. (2006).
\newblock Minimization of {$\phi$}-divergences on sets of signed measures.
\newblock {\em Studia Sci. Math. Hungar.}, {\bf 43}(4), 403--442.

\bibitem[Broniatowski and Keziou(2009)]{BroniatowskiKeziou2009}
Broniatowski, M. and Keziou, A. (2009).
\newblock Parametric estimation and tests through divergences and the duality
  technique.
\newblock {\em J. Multivariate Anal.}, {\bf 100}(1), 16--36.


\bibitem[Broniatowski and Vajda(2009)]{BroniatowskiVajda2009}
Broniatowski, M. and Vajda, I. (2009).
\newblock Several applications of divergence criteria in continuous families.
\newblock Technical Report 2257, Academy of Sciences of the Czech Republic,
  Institute of Information Theory and Automation.



\bibitem[Chatterjee and Bose(2005)]{ChatterjeeBose2005}
Chatterjee, S. and Bose, A. (2005).
\newblock Generalized bootstrap for estimating equations.
\newblock {\em Ann. Statist.}, {\bf 33}(1), 414--436.

\bibitem[Cheng and Huang(2010)]{ChengHuang2010}
Cheng, G. and Huang, J. (2010).
\newblock Bootstrap consistency for general semiparametric
              {$M$}-estimation.
\newblock {\em Ann. Statist.}, {\bf 38}(5), 2884--2915.

\bibitem[Cherfi(2011a)]{Cherfi2011b}
Cherfi, M. (2011a).
\newblock Dual divergences estimation for censored survival data.
\newblock {\em Arxiv preprint: arXiv:1106.2627}.

\bibitem[Cherfi(2011b)]{Cherfi2011a}
Cherfi, M. (2011b).
\newblock Dual $\phi$-divergences estimation in normal models.
\newblock {\em Arxiv preprint: arXiv:1108.2999}.


\bibitem[Cressie and Read(1984)]{CressieRead1984}
Cressie, N. and Read, T. R.~C. (1984).
\newblock Multinomial goodness-of-fit tests.
\newblock {\em J. Roy. Statist. Soc. Ser. B}, {\bf 46}(3), 440--464.

\bibitem[Cs{\"o}rg{\H{o}} and Mason(1989)]{CsorgoMason1989}
Cs{\"o}rg{\H{o}}, S. and Mason, D.~M. (1989).
\newblock Bootstrapping empirical functions.
\newblock {\em Ann. Statist.}, {\bf 17}(4), 1447--1471.

\bibitem[del Barrio and Matr{\'a}n(2000)]{delBarrioMatran2000}
del Barrio, E. and Matr{\'a}n, C. (2000).
\newblock The weighted bootstrap mean for heavy-tailed distributions.
\newblock {\em J. Theoret. Probab.}, {\bf 13}(2), 547--569.

\bibitem[Efron(1979)]{Efron79}
Efron, B. (1979).
\newblock Bootstrap methods: another look at the jackknife.
\newblock {\em Ann. Statist.}, {\bf 7}(1), 1--26.

\bibitem[Efron and Tibshirani(1993)Efron and Tibshirani]{EfronTibshirani1993}
Efron, B. and Tibshirani, R.~J. (1993).
\newblock {\em An introduction to the bootstrap\/}, volume~57 of {\em
  Monographs on Statistics and Applied Probability\/}.
\newblock Chapman and Hall, New York.


\bibitem[G{\"a}nssler(1992)]{Gaenssler1992}
G{\"a}nssler, P. (1992).
\newblock Confidence bands for probability distributions on
  {V}apnik-{C}hervonenkis classes of sets in arbitrary sample spaces using the
  bootstrap.
\newblock In {\em Bootstrapping and related techniques ({T}rier, 1990)}, volume
  376 of {\em Lecture Notes in Econom. and Math. Systems}, pages 57--61.
  Springer, Berlin.

\bibitem[Gin{\'e} and Zinn(1989)]{GineZinn1989}
Gin{\'e}, E. and Zinn, J. (1989).
\newblock Necessary conditions for the bootstrap of the mean.
\newblock {\em Ann. Statist.}, {\bf 17}(2), 684--691.

\bibitem[Gin{\'e} and Zinn(1990)]{GineZinn1990}
Gin{\'e}, E. and Zinn, J. (1990).
\newblock Bootstrapping general empirical measures.
\newblock {\em Ann. Probab.}, {\bf 18}(2), 851--869.

\bibitem[Hong and Kim(2001)]{Hong2001}
Hong, C. and Kim, Y. (2001).
\newblock Automatic selection of the tuning parameter in the minimum density
  power divergence estimation.
\newblock {\em J. Korean Statist. Soc.}, {\bf 30}(3), 453--465.

\bibitem[Holmes and Reinert(2004)]{HolmesReinert2004}
Holmes, S. and Reinert, G. (2004).
\newblock Stein's method for the bootstrap.
\newblock In {\em Stein's method: expository lectures and applications},
  volume~46 of {\em IMS Lecture Notes Monogr. Ser.}, pages 95--136. Inst. Math.
  Statist., Beachwood, OH.

\bibitem[Jim{\'e}nez and Shao(2001)Jim{\'e}nez and Shao]{JimenezShao2001}
Jim{\'e}nez, R. and Shao, Y. (2001).
\newblock On robustness and efficiency of minimum divergence estimators.
\newblock {\em Test\/}, {\bf 10}(2), 241--248.

\bibitem[Keziou(2003)]{Keziou2003}
Keziou, A. (2003).
\newblock Dual representation of {$\phi$}-divergences and applications.
\newblock {\em C. R. Math. Acad. Sci. Paris}, {\bf 336}(10), 857--862.

\bibitem[Kosorok(2008)]{Kosorok2008}
Kosorok, M.~R. (2008).
\newblock {\em Introduction to empirical processes and semiparametric
  inference}.
\newblock Springer Series in Statistics. Springer-Verlag, New York.
\newblock With applications to statistics.


\bibitem[Liese and Vajda(1987)]{LieseVajda1987}
Liese, F. and Vajda, I. (1987).
\newblock {\em Convex statistical distances}, volume~95 of {\em Teubner-Texte
  zur Mathematik [Teubner Texts in Mathematics]}.
\newblock BSB B. G. Teubner Verlagsgesellschaft, Leipzig.
\newblock With German, French and Russian summaries.

\bibitem[Liese and Vajda(2006)]{LieseVajda2006}
Liese, F. and Vajda, I. (2006).
\newblock On divergences and informations in statistics and information theory.
\newblock {\em IEEE Trans. Inform. Theory}, {\bf 52}(10), 4394--4412.

\bibitem[Lo(1993)]{Lo1993}
Lo, A.~Y. (1993).
\newblock A {B}ayesian bootstrap for censored data.
\newblock {\em Ann. Statist.}, {\bf 21}(1), 100--123.

\bibitem[Lohse(1987)]{Lohse1987}
Lohse, K. (1987).
\newblock Consistency of the bootstrap.
\newblock {\em Statist. Decisions}, {\bf 5}(3-4), 353--366.

\bibitem[James(1997)]{Lancelot1997}
James, L.~F. (1997).
\newblock A study of a class of weighted bootstraps for censored data.
\newblock {\em Ann. Statist.}, {\bf 25}(4), 1595--1621.

\bibitem[Kaplan and Meier(1958)]{kaplanmeier1958}
Kaplan, E.~L. and Meier, P. (1958).
\newblock Nonparametric estimation from incomplete observations.
\newblock {\em J. Amer. Statist. Assoc.}, {\bf 53}, 457--481.

\bibitem[Mason and Newton(1992)]{MasonNewton1992}
Mason, D.~M. and Newton, M.~A. (1992).
\newblock A rank statistics approach to the consistency of a general bootstrap.
\newblock {\em Ann. Statist.}, {\bf 20}(3), 1611--1624.

\bibitem[Oakes(1986)]{Oakes1986}
Oakes, D. (1986).
\newblock An approximate likelihood procedure for censored data.
\newblock {\em Biometrics}, {\bf 42}(1), 177--182.

\bibitem[Pardo(2006)]{Pardo2006}
Pardo, L. (2006).
\newblock {\em Statistical inference based on divergence measures}, volume 185
  of {\em Statistics: Textbooks and Monographs}.
\newblock Chapman \& Hall/CRC, Boca Raton, FL.

\bibitem[Pr{\ae}stgaard and Wellner(1993)]{PraestgaardWellner1993}
Pr{\ae}stgaard, J. and Wellner, J.~A. (1993).
\newblock Exchangeably weighted bootstraps of the general empirical process.
\newblock {\em Ann. Probab.}, {\bf 21}(4), 2053--2086.

\bibitem[R{\'e}nyi(1961)]{Reny1961}
R{\'e}nyi, A. (1961).
\newblock On measures of entropy and information.
\newblock In {\em Proc. 4th {B}erkeley {S}ympos. {M}ath. {S}tatist. and
  {P}rob., {V}ol. {I}}, pages 547--561. Univ. California Press, Berkeley,
  Calif.

\bibitem[Rubin(1981)]{Rubin1981}
Rubin, D.~B. (1981).
\newblock The {B}ayesian bootstrap.
\newblock {\em Ann. Statist.}, {\bf 9}(1), 130--134.

\bibitem[Stute(1995)]{Stute1995}
Stute, W. (1995).
\newblock The statistical analysis of {K}aplan-{M}eier integrals.
\newblock In {\em Analysis of censored data ({P}une, 1994/1995)}, volume~27 of
  {\em IMS Lecture Notes Monogr. Ser.}, pages 231--254. Inst. Math. Statist.,
  Hayward, CA.

\bibitem[Toma and Broniatowski(2010)Toma and
  Broniatowski]{TomaBroniatowski2010}
Toma, A. and Broniatowski, M. (2010).
\newblock {Dual divergence estimators and tests: robustness results}.
\newblock {\em J. Multivariate Anal.}, {\bf 102}(1), 20--36.


\bibitem[Tsukahara(2005)]{tsukahara2005}
Tsukahara, H. (2005).
\newblock Semiparametric estimation in copula models.
\newblock {\em Canad. J. Statist.}, {\bf 33}(3), 357--375.


\bibitem[van~der Vaart(1998)]{vanderVaart1998}
van~der Vaart, A.~W. (1998).
\newblock {\em Asymptotic statistics}, volume~3 of {\em Cambridge Series in
  Statistical and Probabilistic Mathematics}.
\newblock Cambridge University Press, Cambridge.

\bibitem[van~der Vaart and Wellner(1996)]{vanderVaartWellner1996}
van~der Vaart, A.~W. and Wellner, J.~A. (1996).
\newblock {\em Weak convergence and empirical processes}.
\newblock Springer Series in Statistics. Springer-Verlag, New York.
\newblock With applications to statistics.

\bibitem[Warwick and Jones(2005)]{Jones2005}
Warwick, J. and Jones, M.~C. (2005).
\newblock Choosing a robustness tuning parameter.
\newblock {\em J. Stat. Comput. Simul.}, {\bf 75}(7), 581--588.

\bibitem[Wellner and Zhan(1996)]{WellnerZhan1996}
Wellner, J,~A. and Zhan, Y. (1996). 
\newblock Bootstrapping $Z$-estimators.
\newblock {\em Technical Report 308, July 1996.}


\bibitem[Weng(1989)]{Weng1989}
Weng, C.-S. (1989).
\newblock On a second-order asymptotic property of the {B}ayesian bootstrap
  mean.
\newblock {\em Ann. Statist.}, {\bf 17}(2), 705--710.

\bibitem[Zheng and Tu(1988)]{Zheng1988}
Zheng, Z.~G. and Tu, D.~S. (1988).
\newblock Random weighting methods in regression models.
\newblock {\em Sci. Sinica Ser. A}, {\bf 31}(12), 1442--1459.


\end{thebibliography}
\end{document}